\documentclass[aos,MSNbibl,seceqn,nameyear,rotating,dvips]{arximspdf}
\usepackage{dcolumn}
\usepackage{multirow}
\usepackage{graphicx}

\doi{10.1214/11-AOS937}
\volume{39}
\issue{6}
\pubyear{2011}
\firstpage{3234}
\lastpage{3261}

\makeatletter

\newcolumntype{d}[1]{D{.}{.}{#1}}

\newtheorem{theorem}{Theorem}[section]
\newtheorem{corollary}[theorem]{Corollary}

\makeatother

\begin{document}
\begin{frontmatter}

\title{Asymptotic properties of the sequential empirical ROC,
PPV and NPV curves under case-control~sampling\thanksref{T1}}
\runtitle{Sequential empirical ROC, PPV and NPV curves}

\thankstext{T1}{Supported by NIH Grants P01-CA53996 and
U01-CA86368.}

\begin{aug}
\author[A]{\fnms{Joseph S.} \snm{Koopmeiners}\corref{}\ead[label=e1]{koopm007@umn.edu}}
\and
\author[B]{\fnms{Ziding} \snm{Feng}\ead[label=e2]{zfeng@fhcrc.org}}
\runauthor{J. S. Koopmeiners and Z. Feng}
\affiliation{University of Minnesota and Fred Hutchinson Cancer
Research Center}
\address[A]{Division of Biostatistics \\
School of Publich Health \\
University of Minnesota \\
A460 Mayo Building, MMC 303 \\
420 Delaware St. SE \\
Minneapolis, Minnesota 55455 \\
USA\\
\printead{e1}}
\address[B]{Division of Public Health Sciences\\
Fred Hutchinson Cancer Research Center \\
M2-B500 \\
PO Box 19024 \\
Seattle, Washington 98109 \\
USA\\
\printead{e2}} %adresu isvedimo komanda gale!
\end{aug}

% HISTORY:
\received{\smonth{3} \syear{2011}}
\revised{\smonth{9} \syear{2011}}

% ABSTRACT
%
\begin{abstract}
The receiver operating characteristic (ROC) curve, the positive
predictive value (PPV) curve and the negative predictive value (NPV)
curve are three measures of performance for a continuous diagnostic
biomarker. The ROC, PPV and NPV curves are often estimated empirically
to avoid assumptions about the distributional form of the biomarkers.
Recently, there has been a push to incorporate group sequential methods
into the design of diagnostic biomarker studies. A~thorough
understanding of the asymptotic properties of the sequential empirical
ROC, PPV and NPV curves will provide more flexibility when designing
group sequential diagnostic biomarker studies. In this paper, we derive
asymptotic theory for the sequential empirical ROC, PPV and NPV curves
under case-control sampling using sequential empirical process theory.
We show that the sequential empirical ROC, PPV and NPV curves converge
to the sum of independent Kiefer processes and show how these results
can be used to derive asymptotic results for summaries of the
sequential empirical ROC, PPV and NPV curves.
\end{abstract}

% KEYWORDS
%

\end{frontmatter}

%s1 #&#
\section{Introduction}
%%%%%%%
Several recent papers have discussed the application of group
sequential methodology to diagnostic biomarker studies
[\citet{TEandZ08}, \citet{TandA10}, \citet{PFLandK09}].
Group sequential study designs (i.e., study
designs with multiple interim analyses) provide an opportunity to
improve the efficiency of diagnostic biomarker studies by allowing
studies to terminate early when the candidate marker is clearly
superior or inferior to established markers or historical levels of
marker performance. Many group sequential methods assume the existence
of a test statistic\vadjust{\goodbreak} with an independent increments covariance structure
[\citet{JandT00}]. A thorough understanding of the asymptotic properties
of the sequential empirical ROC, PPV and NPV curves and, specifically,
verifying that their summary measures have an independent increments
covariance structure, would provide great flexibility when designing
group sequential diagnostic biomarker studies.

Diagnostic biomarkers are used to classify a patient as a case or a
control. A~dichotomous biomarker results in either a positive test,
indicating that the subject should be classified as a case, or a
negative test, indicating that the subject should be classified as a
control. Many biomarkers are measured on a continuous scale and a
threshold must be defined in order to translate a continuous biomarker
into a positive or negative test result. Let $D$ be a~Bernoulli random
variable indicating disease status with prevalence~$\rho$ and let~$X$
be a biomarker value with conditional distribution $F(x \vert D = 1)
\equiv F_{D}(x)$ and $F(x \vert D = 0) \equiv F_{\bar{D}}(x)$, where
$F_{D}(x)$ is the distribution function for the cases and $F_{\bar
{D}}(x)$ is the distribution function for the controls. Furthermore, we
define $F(x) \equiv F_{D}(x) + (1 - \rho) F_{\bar{D}}(x)$
to be the biomarker distribution function for the entire population.
Without loss of generality, assume that larger biomarker values are
more indicative of disease. For a threshold~$c$, a~biomarker value $X$
is translated into a positive test result if it is greater than $c$ and
a negative test result if it is less than or equal to $c$.

The receiver operating characteristic (ROC) curve summarizes the
classification accuracy of a continuous diagnostic biomarker
[\citet{Pepe03}] by reporting the true positive fraction (TPF) and the false
positive fraction (FPF) for all possible cut-offs of the marker. For a
threshold $c$, $\mathrm{TPF}(c) = P [ X > c \vert D = 1 ]$ and
$\mathrm{FPF}(c) = P [ X > c \vert D = 0 ]$. The ROC curve is defined~as
\[
\mathrm{ROC}(c) = \{ (\mathrm{TPF}(c), \mathrm{FPF}(c) ), c \in( -\infty, \infty
) \}
\]
and can alternately be expressed as
%
%e1.1 #&#
%
\begin{equation}
\label{ROCalt}
\mathrm{ROC}(t) = S_{D} (S^{-1}_{\bar{D}} (t ) ),\qquad
t \in(0, 1 ),
\end{equation}
where $S_{D}(x) = 1 - F_{D}(x)$ and $S_{\bar{D}}(x) = 1 - F_{\bar
{D}}(x)$. $\mathrm{ROC}(t)$ can be interpreted as the TPF corresponding to a FPF
of $t$. Alternately, one might be interested in the inverse of the ROC curve,
%
%e1.2 #&#
%
\begin{equation}
\label{ROCinvalt}
\mathrm{ROC}^{-1}(v) = S_{\bar{D}} (S^{-1}_{D} (v ) ),\qquad
v \in(0, 1 ).
\end{equation}
$\mathrm{ROC}^{-1}(v)$ is indexed by the TPF and can be interpreted as the FPF
corresponding to a TPF of $v$.

The predictive accuracy of a dichotomous biomarker can be summarized by
the positive predictive value (PPV) and negative predictive value
(NPV). The PPV and NPV curves were proposed as an extension of PPV and
NPV to continuous markers [\citet{MandP04},
\citet{ZCPandL08}]. For a threshold $c$, $\mathrm{PPV}(c) = P [ D = 1 \vert
X > c ]$ and $\mathrm{NPV}(c) = P [ D = 0 \vert X \leq c ]$. The PPV and NPV
curves are defined as $\mathrm{PPV}(c)$ and $\mathrm{NPV}(c)$ for all $c \in( - \infty,
\infty)$. In practice, PPV and NPV curves are indexed by a summary of
the marker distribution rather than a generic threshold
[\citet{MandP04}, \citet{ZCPandL08}]. In this paper, we
consider the PPV and NPV curves indexed by the FPF and the percentile
value in the entire population.

The ROC, PPV and NPV curves are commonly estimated nonparametrically to
avoid making assumptions about the form of $F_{D}(x)$ and $F_{\bar
{D}}(x)$. This is particularly important in the case of the ROC, PPV
and NPV curves because we are often interested in regions of the curve
that correspond to the tails of these distributions. For example, a
biomarker must possess a~high specificity in order to be clinically
useful in a low disease risk population screening setting, which
corresponds to the upper tail of the biomarker distribution among
controls.

Our understanding of the empirical ROC curve is enhanced by knowledge
of its asymptotic properties. \citet{HandT96} showed that the empirical
ROC curve converges to the sum of two independent Brownian bridges. The
asymptotic normality of summary measures of the empirical ROC curve,
such as the area under the ROC curve or a point on the ROC curve, can
be derived from their work. To our knowledge, no asymptotic theory is
available for the empirical PPV and NPV curves.

\citet{TEandZ08} showed that a family of weighted area under the ROC
curve (wAUC) statistics has an independent increments covariance
structure. It would be beneficial to show that this assumption holds
for a larger class of summaries of the ROC curve. In this paper, we
develop asymptotic theory for the sequential empirical ROC, PPV and NPV
curves. Our results allow us to develop distribution theory for other
summaries of the ROC curve and to develop distribution theory for
summaries of the PPV and NPV curves.

%s2 #&#
\section{Notation and definitions}
Before beginning our discussion of the sequential empirical ROC, PPV
and NPV curves, we provide definitions of the sequential empirical
estimates for the underlying distribution and quantile functions. Let
$X_{D,1}, X_{D,2}, \ldots, X_{D,n_{D}}$ be i.i.d. marker values for
the cases with distribution function, $F_{D}(x)$, and $X_{\bar{D},1},
X_{\bar{D},2}, \ldots, X_{\bar{D},n_{\bar{D}}}$ be i.i.d. marker
values for the controls with distribution function, $F_{\bar{D}}(x)$.
Furthermore, let $r_{D}$ and $r_{\bar{D}}$ refer to the proportion of
case and controls, respectively, that are observed at a given time
point. The sequential empirical estimate of $F_{D}(x)$ is defined as
\[
\hat{F}_{D, r_{D}}(x) = \cases{
0, &\quad$\displaystyle 0 \leq r_{D} < \frac{1} {n_{D}}$,\vspace*{1pt}\cr
\displaystyle \frac{1}{[r_{D}n_{D}]} \sum_{i = 1}^{[r_{D} n_{D}]} 1\{X_{D,i} \leq x
\}, &\quad$\displaystyle - \infty< x < \infty, \frac{1}{n_{D}} \leq r_{D}
\leq1$,}\vadjust{\goodbreak}
\]
and the sequential empirical estimate of $F_{D}^{-1}(t)$ is defined as
\[
\hat{F}_{D, r_{D}}^{-1} (t) = \cases{
X_{D,1,[r_{D}n_{D}]}, &\quad if $t = 0, 0 \leq
r_{D} \leq1$, \vspace*{2pt}\cr
X_{D,k,[r_{D}n_{D}]}, &\quad if $\displaystyle \frac{k-1}{[r_{D}n_{D}]} < t
\leq\frac{k}{[r_{D}n_{D}]}$, \vspace*{2pt}\cr
&\quad$1 \leq k \leq[r_{D}n_{D}], 0 \leq t \leq1$,}
\]
where $X_{D,1,[r_{D}n_{D}]}, X_{D,2,[r_{D}n_{D}]}, \ldots,
X_{D,[r_{D}n_{D}],[r_{D}n_{D}]}$ are the sequential order statistics of
the biomarker values for the cases. The sequential empirical estimates
of $S_{D}(x)$ and $S_{D}^{-1}(t)$ are defined as $\hat{S}_{D,r_{D}}(x)
= 1 - \hat{F}_{D, r_{D}}(x)$ and $\hat{S}_{D,r_{D}}^{-1} (t) = \hat
{F}_{D, r_{D}}^{-1} (1-t)$. The sequential empirical estimates for the
control population are defined in an analogous fashion. The sequential
empirical estimates of $F_{D}(x)$ and $F_{\bar{D}}(x)$ lead to a
natural definition of the sequential empirical estimates of $F(x)$ and
$F^{-1}(t)$,
\[
\hat{F}_{r_{D}, r_{\bar{D}}}(x) = \rho\hat{F}_{D, r_{D}}(x) + (1
- \rho) \hat{F}_{\bar{D}, r_{\bar{D}}}(x)
\]
and
\[
\hat{F}_{r_{D}, r_{\bar{D}}}^{-1}(t) = \inf\{x\dvtx\hat{F}_{r_{D},
r_{\bar
{D}}}(x) \geq t \},
\]
where $\rho$ is assumed to be known. $\hat{F}_{r_{D}, r_{\bar{D}}}(x)$
is a linear combination of~$\hat{F}_{D, r_{D}}(x)$ and $\hat{F}_{\bar
{D}, r_{\bar{D}}}(x)$ and is therefore indexed by both $r_{D}$, the
proportion of cases observed at a given time point, and $r_{\bar{D}}$,
the proportion of controls observed at a given time point.

Throughout this paper, we let $0 < a < b < 1$, $0 < c < 1$, $0 < d < 1$
and make the following assumptions:
\begin{longlist}[(A4)]
\item[(A1)] $F_{D}(x)$ and $F_{\bar{D}}(x)$ are continuous
distribution functions with continuous densities $f_{D}(x)$ and $f_{\bar
{D}}(x)$, respectively,
\item[(A2)] $f_{D}(x) > 0 $ for $x \in(\sup\{x\dvtx F_{D}(x) = 0\}
, \inf\{x\dvtx F_{D}(x) = 1\} )$,
\item[(A3)] $f_{\bar{D}}(x) > 0 $ for $x \in(\sup\{x\dvtx F_{\bar
{D}}(x) = 0\}, \inf\{x\dvtx F_{\bar{D}}(x) = 1\} )$,\vspace*{1pt}
\item[(A4)] $\frac{n_{D}} {n_{\bar{D}}}\rightarrow\lambda> 0$ as
$n_{D} \rightarrow\infty$ and $n_{\bar{D}} \rightarrow\infty$, that
is, the ratio of cases to controls converges to a constant that is
greater than 0.
\end{longlist}
The asymptotic results in Section \ref{asympsec} make use of the
Kiefer process. The Kiefer process, $K(t,r)$, is a two-dimensional,
mean-zero Gaussian process with covariance
\[
\operatorname{Cov}(K(t_{1}, r_{1}), K(t_{2}, r_{2}))
= ( t_{1} \wedge t_{2} - t_{1} t_{2} ) ( r_{1} \wedge
r_{2} ),
\]
where $\wedge$ represents the minimum. The Kiefer process behaves like
a Brownian bridge in $t$ and Brownian Motion in $r$.

The remainder of this paper proceeds as follows. In Section \ref
{asympsec}, we develop asymptotic theory for the sequential empirical
ROC, PPV and NPV curves. First, we generalize the work of \citet
{HandT96} to the sequential empirical ROC curve by showing that the
sequential empirical ROC curve converges to the sum of independent
Kiefer processes. Next, we develop asymptotic theory for the sequential
empirical\vadjust{\goodbreak} PPV and NPV curves indexed by the FPF by writing them as
functions of the sequential empirical ROC curve. Finally, we follow the
approach of \citet{PandS68} to develop asymptotic theory for the
PPV and
NPV curves indexed by the percentile value of the marker distribution.
We validate our asymptotic results by simulation in Section \ref
{simsec} and illustrate how they can be used to design group
sequential diagnostic biomarker studies in Section \ref{appsec}. We
conclude with a discussion in Section \ref{discsec}.\vspace*{-2pt}

%s3 #&#
\section{Asymptotic results}\vspace*{-2pt}
\label{asympsec}

%s3.1 #&#
\subsection{The sequential empirical ROC curve}
\label{ROCsec}
%%%%%%%
In this section, we provide asymptotic results for the sequential
empirical ROC curve. Results for the inverse of the sequential
empirical ROC curve are nearly identical; we direct the reader to an
associated technical report for details [\citet{KandF10}]. The sequential
empirical ROC curve, $\widehat{\mathrm{ROC}}_{r_{D}, r_{D}}(t)$, is defined by
substituting the sequential empirical estimates of $S_{D}(x)$ and
$S_{\bar{D}}(x)$ into~(\ref{ROCalt}), yielding
\[
\widehat{\mathrm{ROC}}_{r_{D},r_{\bar{D}}}(t) = \hat{S}_{D, r_{D}} (\hat
{S}^{-1}_{\bar{D},r_{\bar{D}}} (t ) ),
\]
and for ease of notation, we define
\[
R_{r_{D}, r_{\bar{D}}}(t) \equiv n_{D}^{-1/2} [n_{D} r_{D}] \bigl(\widehat
{\mathrm{ROC}}_{r_{D}, r_{\bar{D}}}(t) - \mathrm{ROC}(t)\bigr).
\]
The primary result in this section provides asymptotic theory for
$R_{r_{D}, r_{\bar{D}}}(t)$. By developing asymptotic theory for
$R_{r_{D}, r_{\bar{D}}}(t)$, we are also able to develop asymptotic
theory for functionals of $R_{r_{D}, r_{\bar{D}}}(t)$ as a special
case. Theorem \ref{ROCseqproc} establishes the convergence of
$R_{r_{D}, r_{\bar{D}}}(t)$ to the sum of independent Kiefer processes.
%
%th3.1 #&#
%
\begin{theorem}
\label{ROCseqproc}
%%%%%%%
Assume \textup{(A1)--(A4)} hold and let $\frac{f_{D}(S_{\bar{D}}^{-1}(t))} {f_{\bar
{D}}(S_{\bar{D}}^{-1}(t)) }$ be bounded on $[a,b]$. As $n_{D}
\rightarrow\infty$ and $n_{\bar{D}} \rightarrow\infty$
\[
R_{r_{D}, r_{\bar{D}}}(t) \rightarrow_{d} K_1(\mathrm{ROC}(t), r_{D}) +
\lambda^{1/2} \frac{r_{D}}{r_{\bar{D}}}
\biggl(\frac{f_{D}(S_{\bar{D}}^{-1}(t))}
{f_{\bar{D}}(S_{\bar{D}}^{-1}(t)) } \biggr) K_2(t, r_{\bar{D}})
\]
uniformly for $t \in[a, b ]$, $r_{D} \in[c, 1
]$ and $r_{\bar{D}} \in[d, 1 ]$ where $K_{1}$ and $K_{2}$
are independent Kiefer processes.
%%%%%%%
\end{theorem}

%%%%%%%
A proof of Theorem \ref{ROCseqproc} can be found in the \hyperref[app]{Appendix}.
Theorem \ref{ROCseqproc} generalizes the results of \citet
{HandT96} to
the sequential empirical ROC curve. The proof of Theorem \ref
{ROCseqproc} is similar to the proof found in \citet{HandT96} but our
proof relies on the more powerful sequential empirical process theory.
Sequential empirical process theory generalizes asymptotic theory for
the standard empirical process by introducing a parameter for time. In
doing so, asymptotic results for the sequential empirical process\vadjust{\goodbreak}
involve the Kiefer process. Using properties of the Kiefer process, we
are able to easily derive asymptotic results for summaries of the
sequential empirical ROC curve and verify that the independent
increments assumption holds in many cases. Furthermore,
we can recover Hsieh and Turnbull's result as a special case of
Theorem \ref{ROCseqproc} by letting $r_{D}$ and $r_{\bar{D}}$ both
equal~1.

%co3.2 #&#
%
\begin{corollary}
%%%%%%%
\label{ROCfixedproc}
%%%%%%%
Assume \textup{(A1)--(A4)} hold and let $\frac{f_{D}(S_{\bar{D}}^{-1}(t))} {f_{\bar
{D}}(S_{\bar{D}}^{-1}(t)) }$ be bounded on $[a,b]$. As $n_{D}
\rightarrow\infty$ and $n_{\bar{D}} \rightarrow\infty$,
\[
%n_{D}^{1/2} (\widehat{\mathrm{ROC}}_{1,1}(t) - \mathrm{ROC}(t)) \rightarrow_{d} &
%B_1(\mathrm{ROC}(t)) + \lambda^{1/2} (\frac{f_{D}(S_{\bar{D}}^{-1}(t))}
%{f_{\bar{D}}(S_{\bar{D}}^{-1}(t)) } ) B_2(t)
R_{1,1}(t) \rightarrow_{d} B_1(\mathrm{ROC}(t)) + \lambda^{1/2}
\biggl(\frac{f_{D}(S_{\bar{D}}^{-1}(t))}
{f_{\bar{D}}(S_{\bar{D}}^{-1}(t)) } \biggr) B_2(t)
\]
uniformly for $t \in[a, b ]$ where $B_{1}$ and $B_{2}$ are
independent Brownian bridges.
%%%%%%%
\end{corollary}
\begin{pf}
%%%%%%%
Immediate from Theorem \ref{ROCseqproc} and by noting that $K(t, 1)
=_{d} B(t)$.
%%%%%%%
\end{pf}

An advantage to studying the asymptotic behavior of the sequential
empirical ROC curve at the process level, rather than a single point on
the sequential empirical ROC curve, is that we are able to study the
joint behavior of multiple points on the ROC curve. Corollary \ref
{ROCmultt} provides a normal approximation for a vector of points on
the sequential empirical ROC curve.

%co3.3 #&#
%
\begin{corollary}
%%%%%%%
\label{ROCmultt}
%%%%%%%
Assume \textup{(A1)--(A4)} hold and let
$\frac{f_{D}(S_{\bar{D}}^{-1}(t))} {f_{\bar {D}}(S_{\bar{D}}^{-1}(t))
}$ be bounded on $[a,b]$. For\vspace*{1pt} $t_{1}, t_{2}, \ldots, t_{J}
\in(0,1)$, $r_{D,1}, r_{D,2}, \ldots, r_{D,J} \in(0,1]$ and
$r_{\bar{D},1}, r_{\bar{D},2},\allowbreak \ldots, r_{\bar {D},J} \in(0,1]$,
a vector of arbitrary points on the sequential empirical ROC curve,
$(\widehat{\mathrm{ROC}}_{r_{D,1}, r_{\bar {D},1}}(t_{1}),
\widehat{\mathrm{ROC}}_{r_{D,2}, r_{\bar{D},2}}(t_{2}), \ldots,
\widehat{\mathrm{ROC}}_{r_{D,J}, r_{\bar{D},J}}(t_{J}) )$, is
approximately multivariate normal with
\[
\widehat{\mathrm{ROC}}_{r_{D,j}, r_{\bar{D},j}}(t_{j}) \sim N \bigl(\mathrm{ROC}(t_{j}),
\sigma^{2}_{\widehat{\mathrm{ROC}}_{r_{D,j}, r_{\bar{D},j}}(t_{j})}
\bigr),\qquad
j = 1, 2, \ldots, J,
\]
where
\[
\sigma^{2}_{\widehat{\mathrm{ROC}}_{r_{D,j}, r_{\bar{D},j}}(t_{j})} = \frac{
\mathrm{ROC}(t_{j}) (1 - \mathrm{ROC}(t_{j}) ) } {n_{D} r_{D,j} } +
\biggl(\frac{f_{D}(S_{\bar{D}}^{-1}(t_{j}))}
{f_{\bar{D}}(S_{\bar{D}}^{-1}(t_{j})) } \biggr)^2 \frac{t_{j} (1 -
t_{j} ) } {n_{\bar{D}} r_{\bar{D},j}}
\]
and
\begin{eqnarray*}
&&\operatorname{Cov} [\widehat{\mathrm{ROC}}_{r_{D,i}, r_{\bar{D},i}}(t_{i}),
\widehat{\mathrm{ROC}}_{r_{D,j}, r_{\bar{D},j}}(t_{j}) ] \\[-2pt]
&&\qquad = \frac{ (r_{D,i} \wedge r_{D,j} ) ( \mathrm{ROC}(t_{i})
\wedge \mathrm{ROC}(t_{j}) - \mathrm{ROC}(t_{i}) \operatorname{ROC}(t_{j}) ) } {n_{D} r_{D,i}
r_{D,j} } \\[-2pt]
&&\qquad\quad{} + \biggl(\frac{f_{D}(S_{\bar{D}}^{-1}(t_{i}))} {f_{\bar{D}}(S_{\bar
{D}}^{-1}(t_{i})) } \biggr) \biggl(\frac{f_{D}(S_{\bar{D}}^{-1}(t_{j}))}
{f_{\bar{D}}(S_{\bar{D}}^{-1}(t_{j})) } \biggr) \frac{ ( r_{\bar
{D},i} \wedge r_{\bar{D},j} ) ( t_{i} \wedge t_{j} - t_{i}
t_{j} ) } {n_{\bar{D}} r_{\bar{D},i} r_{\bar{D},j}}.\vspace*{-3pt}
\end{eqnarray*}
\end{corollary}
\begin{pf}
%%%%%%%
Immediate from Theorem \ref{ROCseqproc}.\vadjust{\goodbreak}
%%%%%%%
\end{pf}

Corollary \ref{ROCmultt} provides the asymptotic covariance for two
points at different locations and different times on the sequential
empirical ROC curve. This allows us to fully specificy the joint
sequential distribution of multiple points on the ROC curve, which
allows us to design group sequential diagnostic biomarker studies where
multiple points on the ROC curve are treated as multiple endpoints of a
group sequential study. For example, we might be interested in
$\mathrm{ROC}(t_{1})$ and $\mathrm{ROC}(t_{2})$, where $t_{1}$ is chosen for high
specificity to rule patients in for work up and $t_{2}$ is chosen for
high sensitivity to rule out patients for invasive work.

Our interest in the sequential empirical ROC curve is motivated by the
need to design group sequential diagnostic biomarker studies. Our
ability to design group sequential diagnostic biomarker studies would
be enhanced by showing that summaries of the sequential empirical ROC
curve have an independent increments covariance structure. The simplest
summary of the ROC curve is a point on the ROC curve, $\mathrm{ROC}(t)$.
$\mathrm{ROC}(t)$ can be interpreted as the sensitivity at a specificity of $1 -
t$. Corollary \ref{ROCseq} shows that the sequential empirical
estimator of $\mathrm{ROC}(t)$ is asymptotically normal and has independent
increments when divided its variance.\vspace*{-3pt}
%
%co3.4 #&#
%
\begin{corollary}
%%%%%%%
\label{ROCseq}
%%%%%%%
Assume \textup{(A1)--(A4)} hold and let
$\frac{f_{D}(S_{\bar{D}}^{-1}(t))} {f_{\bar {D}}(S_{\bar{D}}^{-1}(t))
}$ be bounded on $[a,b]$. For $t \in(0, 1 )$ and J stopping times,
$(\widehat{\mathrm{ROC}}_{r_{D,1}, r_{\bar {D},1}}(t),
\widehat{\mathrm{ROC}}_{r_{D,2}, r_{\bar{D},2}}(t),\allowbreak \ldots,
\widehat {\mathrm{ROC}}_{r_{D,J}, r_{\bar{D},J}}(t) )$, is
approximately multivariate normal with
\[
\widehat{\mathrm{ROC}}_{r_{D,i}, r_{\bar{D},i}}(t) \sim N \bigl(\mathrm{ROC}(t), \sigma
^{2}_{\widehat{\mathrm{ROC}}_{r_{D,i}, r_{\bar{D},i}}(t)} \bigr),\qquad
i = 1, 2, \ldots, J,\vspace*{-2pt}
\]
and
\begin{eqnarray*}
&&\operatorname{Cov} [\widehat{\mathrm{ROC}}_{r_{D,i}, r_{\bar{D},i}}(t),
\widehat{\mathrm{ROC}}_{r_{D,j}, r_{\bar{D},j}}(t) ] \\[-2pt]
&&\qquad= \operatorname{Var} [ \widehat{\mathrm{ROC}}_{r_{D,j}, r_{\bar{D},j}}(t) ] =
\sigma^{2}_{\widehat{\mathrm{ROC}}_{r_{D,j}, r_{\bar{D},j}}(t)},\qquad
r_{i} \leq r_{j},\vspace*{-2pt}
\end{eqnarray*}
where $\sigma^{2}_{\widehat{\mathrm{ROC}}_{r_{D,j}, r_{\bar{D},j}}(t)}$ is
defined as in Corollary \ref{ROCmultt}.\vspace*{-3pt}
%%%%%%%
\end{corollary}
\begin{pf}
%%%%%%%
Immediate from Corollary \ref{ROCmultt}.\vspace*{-3pt}
%%%%%%%
\end{pf}

Asymptotic theory for other summary measures of the ROC curve, such as
the area under the curve or the partial area under the curve, can also
be derived from Theorem \ref{ROCseqproc}. This illustrates the
flexibility of Theorem \ref{ROCseqproc}. By developing distribution
theory for the sequential empirical ROC curve, we are able to derive
distribution theory for summaries of the ROC curve as a~special case.\vspace*{-3pt}

%s3.2 #&#
\subsection{The sequential empirical PPV and NPV curves indexed by the
false positive fraction}
%%%%%%%
\label{PPVtsec}
%%%%%%%
In this section, we consider the sequential empirical PPV and NPV\vadjust{\goodbreak}
curves indexed by the false positive fraction, $t$. The PPV and NPV
curve indexed by the false positive fraction can be written as a
function of the ROC curve and their asymptotic properties can be
derived using the results from Section \ref{ROCsec}. Asymptotic
results for the PPV and NPV curve indexed by the true positive
fraction, $v$, can similarly be derived by writing the PPV and NPV
curve as a function of the inverse of the ROC curve but are not
presented in this paper. The interested reader is directed to
\citet
{KandF10} for details.

The PPV and NPV curves indexed by the false positive fraction are
defined as $\mathrm{PPV}(t) = P[D = 1 \vert X > S_{\bar{D}}^{-1}(t)]$ and
$\mathrm{NPV}(t) = P[D = 0 \vert X \leq S_{\bar{D}}^{-1}(t)]$ for all $t \in
(0,1)$ and can be written as functions of the ROC curve as follows:
%
%e3.1 #&#
%
\begin{equation}
\label{PPVt}
\mathrm{PPV}(t) %= & P[D = 1 \vert X > S_{\bar{D}}^{-1}(t)] \nonumber\\
%& = \frac{ P[D = 1, X > S_{\bar{D}}^{-1}(t) ] } { P[X > S_{
%& = \frac{ P[X > S_{\bar{D}}^{-1}(t) \vert D = 1] P[D = 1] } { P[X >
%S_{\bar{D}}^{-1}(t) \vert D = 1] P[D = 1] + P[X > S_{\bar{D}}^{-1}(t)
= \frac{ \mathrm{ROC}(t) \rho} {\mathrm{ROC}(t) \rho+ t
(1 - \rho) }
\end{equation}
and
%
%e3.2 #&#
%
\begin{equation}
\label{NPVt}
\mathrm{NPV}(t) = \frac{ (1 - t ) ( 1 - \rho) } { (
1 - \mathrm{ROC}(t) ) \rho+ ( 1 - t ) ( 1-
\rho) }.
\end{equation}
The sequential empirical estimators of $\mathrm{PPV}(t)$ and $\mathrm{NPV}(t)$ are
defined be plugging the sequential empirical estimator of $\mathrm{ROC}(t)$ into
(\ref{PPVt}) and (\ref{NPVt}), yielding
\[
\widehat{\mathrm{PPV}}_{r_{D}, r_{\bar{D}}}(t) = \frac{ \widehat{\mathrm{ROC}}_{r_{D},
r_{\bar{D}}} (t) \rho} {\widehat{\mathrm{ROC}}_{r_{D}, r_{\bar{D}}}
(t) \rho+ t (1 - \rho) }
\]
and
\[
\widehat{\mathrm{NPV}}_{r_{D}, r_{\bar{D}}}(t) = \frac{ (1 - t )
(1 - \rho) } {(1 - \widehat{\mathrm{ROC}}_{r_{D}, r_{\bar{D}}}
(t) ) \rho+ (1 - t ) (1 - \rho) }.
\]
From this point forward, we only consider $\widehat{\mathrm{PPV}}_{r_{D}, r_{\bar
{D}}}(t)$ and note that results for $\widehat{\mathrm{NPV}}_{r_{D}, r_{\bar
{D}}}(t)$ are nearly identical. Again, for ease of notation, we define
\[
P_{r_{D}, r_{\bar{D}}}(t) \equiv n_{D}^{-1/2} [n_{D} r_{D}] \bigl(\widehat
{\mathrm{PPV}}_{r_{D}, r_{\bar{D}}}(t) - \mathrm{PPV}(t)\bigr).
\]
We begin by using the results of Section \ref{ROCsec} to derive
asymptotic theory for~$P_{r_{D},r_{\bar{D}}}(t)$. Theorem \ref
{PPVtseqproc} establishes the convergence of $P_{r_{D},r_{\bar
{D}}}(t)$ to the sum of two independent Kiefer processes.
%
%th3.5 #&#
%
\begin{theorem}
%%%%%%%
\label{PPVtseqproc}
%%%%%%%
Assume \textup{(A1)--(A4)} hold and let $\frac{f_{D}(S_{\bar{D}}^{-1}(t))} {f_{\bar
{D}}(S_{\bar{D}}^{-1}(t)) }$ be bounded on $[a,b]$. As $n_{D}
\rightarrow\infty$ and $n_{\bar{D}} \rightarrow\infty$
\begin{eqnarray*}
P_{r_{D},r_{\bar{D}}}(t) &\rightarrow_{d}& \biggl(\frac{t (1 -
\rho) \rho} { (\mathrm{ROC}(t) \rho+ t ( 1 - \rho)
)^2 }\biggr) \\
&&{}\times\biggl( K_1(\mathrm{ROC}(t), r_{D}) + \lambda^{1/2} \frac
{r_{D}}{r_{\bar{D}}} \biggl(\frac{f_{D}(S_{\bar{D}}^{-1}(t))}
{f_{\bar{D}}(S_{\bar{D}}^{-1}(t)) } \biggr) K_2(t, r_{\bar{D}}) \biggr)
%P_{r_{D},r_{\bar{D}}}(t) \rightarrow_{d} & (\frac{t (1 -
%)^2 }) K_1(\mathrm{ROC}(t), r_{D}) \\
%& + (\frac{t (1 - \rho) \rho} { (\mathrm{ROC}(t) \rho+
%t ( 1 - \rho) )^2 }) \lambda^{1/2} \frac
%{r_{D}}{r_{\bar{D}}} (\frac{f_{D}(S_{\bar{D}}^{-1}(t))} {f_{
\end{eqnarray*}
uniformly for $t \in[a, b ]$, $r_{D} \in[c, 1
]$ and $r_{\bar{D}} \in[d, 1 ]$ where $K_{1}$ and $K_{2}$
are independent Kiefer processes.
%%%%%%%
\end{theorem}

%%%%%%%
The proof of Theorem \ref{PPVtseqproc} relies on writing
$P_{r_{D},r_{\bar{D}}}(t)$ as a function\break of~$R_{r_{D},r_{\bar{D}}}(t)$
\begin{eqnarray*}
P_{r_{D},r_{\bar{D}}}(t) &=&  \biggl( \frac{ \widehat{\mathrm{ROC}}_{r_{D},
r_{\bar{D}}} (t) \rho} {\widehat{\mathrm{ROC}}_{r_{D}, r_{\bar{D}}}
(t) \rho+ t (1 - \rho) } - \frac{ \mathrm{ROC}(t) \rho
} {\mathrm{ROC}(t) \rho+ t (1 - \rho) } \biggr)\\
&&{} \times\bigl({ \widehat
{\mathrm{ROC}}_{r_{D}, r_{\bar{D}}} (t) - \mathrm{ROC}(t) }\bigr)^{-1} R_{r_{D},r_{\bar{D}}}(t)
\end{eqnarray*}
and applying the results of Theorem \ref{ROCseqproc}. The first term
converges to
\[
\biggl(\frac{t (1 - \rho) \rho} {
(\mathrm{ROC}(t) \rho+ t ( 1 - \rho) )^2 }\biggr)
\]
and
$R_{r_{D},r_{\bar{D}}}(t)$ converges to the sum of two independent
Kiefer process by Theorem~\ref{ROCseqproc}. A formal proof of
Theorem \ref{PPVtseqproc} can be found in \citet{KandF10}.

From Theorem \ref{PPVtseqproc}, we can prove analogous results to
Corollaries \ref{ROCmultt} and~\ref{ROCseq} for the sequential
empirical PPV curve indexed by the FPF. Namely, that an arbitrary
vector of points on the sequential empirical PPV curve follows a
multivariate normal distribution and the sequential empirical estimate
of a point on the PPV curve is approximately normally distributed with
an independent increments covariance structure. We leave the formal
statement of these corollaries for the \hyperref[app]{Appendix} but present the form of
the covariance between two arbitrary points on the sequential empirical
PPV curve:\looseness=-1
\begin{eqnarray*}
&& \operatorname{Cov} [\widehat{\mathrm{PPV}}_{r_{D,i},
r_{\bar{D},i}}(t_{i}), \widehat
{\mathrm{PPV}}_{r_{D,j}, r_{\bar{D},j}}(t_{j}) ] \\
&&\qquad = \biggl(\frac{t_{i} (1 - \rho) \rho} { (\mathrm{ROC}(t_{i})
\rho+ t_{i} ( 1 - \rho) )^2 } \biggr) \biggl(\frac
{t_{j} (1 - \rho) \rho} { (\mathrm{ROC}(t_{j}) \rho+ t_{j}
( 1 - \rho) )^2 } \biggr)\\
&&\qquad\quad{}\times \operatorname{Cov} [\widehat
{\mathrm{ROC}}_{r_{D,i}, r_{\bar{D},i}}(t_{i}), \widehat{\mathrm{ROC}}_{r_{D,j}, r_{\bar
{D},j}}(t_{j}) ] .
\end{eqnarray*}\looseness=0
$\mathrm{PPV}(t)$ is a function of $\mathrm{ROC}(t)$ and, therefore, distribution theory
for a~vector of points on the PPV curve can also be derived using the
delta method and Corollary~\ref{ROCmultt}.

Asymptotic theory for the fixed-sample empirical PPV curve indexed by
the FPF, which was previously unavailable, can be derived as a special
case of Theorem \ref{PPVtseqproc} by letting $r_{D}$ and $r_{\bar
{D}}$ equal 1. The fixed-sample empirical PPV curve converges to the
sum of independent Brownian bridges
\begin{eqnarray*}
P_{1,1}(t) &\rightarrow_{d}& \biggl(\frac{t (1 - \rho) \rho}
{ (\mathrm{ROC}(t) \rho+ t ( 1 - \rho) )^2 } \biggr)\\
&&{}\times\biggl( B_1(\mathrm{ROC}(t))+ \lambda^{1/2} \biggl(\frac{f_{D}(S_{\bar
{D}}^{-1}(t))} {f_{\bar{D}}(S_{\bar{D}}^{-1}(t)) } \biggr) B_2(t) \biggr),
\end{eqnarray*}
which allows us to derive a normal approximation for the empirical
estimate of a point on the PPV curve
\[
\widehat{\mathrm{PPV}}_{1,1}(t) \sim N\biggl( \mathrm{PPV}(t), \biggl(\frac{t (1 -
\rho) \rho} { (\mathrm{ROC}(t) \rho+ t ( 1 - \rho)
)^2 } \biggr)^2 \sigma^{2}_{\widehat{\mathrm{ROC}}_{1,1}(t)} \biggr),
\]
where $\sigma^{2}_{\widehat{\mathrm{ROC}}_{1,1}(t)}$ is defined as in
Corollary \ref{ROCmultt}.

%
%s3.3 #&#
\subsection{The sequential empirical PPV and NPV curves indexed by the
percentile value}
%%%%%%%
\label{PPVusec}
%%%%%%%
Finally, we consider the PPV and NPV curves indexed by the proportion
of the population that are classified as negative, $u$, and positive,
$1 - u$. In this case, the PPV and NPV curves are defined as $\mathrm{PPV}(u) =
P[D = 1 \vert X > F^{-1}(u)]$ and $\mathrm{NPV}(u) = P[D = 0 \vert X \leq
F^{-1}(u)]$ for all $u \in(0,1)$. Under this indexing, the PPV curve
can be written as
%
%e3.3 #&#
%
\begin{equation}
\label{PPVu}
\mathrm{PPV}(u) %& = P [ D = 1 \vert X > F^{-1}(u) ] \nonumber\\
%& = \frac{P [ D = 1, X > F^{-1}(u) ] } {P [X >
%F^{-1}(u) ] } \nonumber\\
%& = \frac{P [X > F^{-1}(u) \vert D = 1 ] * P [ D = 1
%]} {1 - F ( F^{-1}(u) )} \nonumber\\
= \frac{S_{D} ( F^{-1}(u) ) \rho} {1 - u},
\end{equation}
and since the NPV curve can be written as
%
%e3.4 #&#
%
\begin{equation}
\label{NPVu}
\mathrm{NPV}(u) = \frac{u - \rho} {u } + \frac{1 - u } {u} \operatorname{PPV}(u),
\end{equation}
it suffices to study the PPV curve when considering estimation of the
NPV curve.

The sequential empirical estimator of $\mathrm{PPV}(u)$ is found by substituting
the sequential empirical estimators of $S_{D}(x)$ and $F(x)$, along
with the known value of~$\rho$, into (\ref{PPVu}),
%
%e3.5 #&#
%
\begin{equation}
\label{PPVuseq}
\widehat{\mathrm{PPV}}_{r_{D}, r_{\bar{D}}}(u) = \frac{\hat{S}_{D, r_{D}}
( \hat{F}_{r_{D}, r_{\bar{D}}}^{-1}(u) ) \rho} {1 - u},
\end{equation}
and the sequential empirical estimator of $\mathrm{NPV}(u)$ is found by
substituting the sequential empirical estimator of $\mathrm{PPV}(u)$ into (\ref{NPVu}),
%
%e3.6 #&#
%
\begin{equation}
\label{NPVuseq}
\widehat{\mathrm{NPV}}_{r_{D}, r_{\bar{D}}}(u) = \frac{u - \rho} {u } + \frac
{1 - u } {u} \,\widehat{\mathrm{PPV}}_{r_{D}, r_{\bar{D}}}(u).
\end{equation}
Finally, we define,
\[
P_{r_{D}, r_{\bar{D}}}(u) = n_{D}^{-1/2} [n_{D} r_{D}] \bigl( \widehat
{\mathrm{PPV}}_{r_{D}, r_{\bar{D}}}(u) - \mathrm{PPV}(u) \bigr)
\]
and
\[
N_{r_{D}, r_{\bar{D}}}(u) = n_{D}^{-1/2} [n_{D} r_{D}] \bigl( \widehat
{\mathrm{NPV}}_{r_{D}, r_{\bar{D}}}(u) - \mathrm{NPV}(u) \bigr)
\]
for mathematical convenience. We begin by developing distribution
theory for $P_{r_{D}, r_{\bar{D}}}(u)$. Theorem \ref{PPVuseqproc}
establishes the convergence of the sequential empirical PPV curve to
the sum of two independent Kiefer processes.
%
%th3.6 #&#
%
\begin{theorem}
%%%%%%%
\label{PPVuseqproc}
%%%%%%%
Assume \textup{(A1)--(A4)} hold and let $\frac{ f_{D} (F^{-1} (u)
) } {f (F^{-1} (u) ) }$ be bounded on $
[a,b]$. As $n_{D} \rightarrow\infty$ and $n_{\bar{D}}
\rightarrow\infty$
\begin{eqnarray*}
P_{r_{D}, r_{\bar{D}}}(u) &\rightarrow_{d}& - \frac{\rho( 1 -
\rho) } {1-u} \frac{ f_{\bar{D}} (F^{-1} (u)
) } {f (F^{-1} (u) ) }
K_1(F_{D}(F^{-1}(u)), r_{D}) \\
&&{} + \frac{\rho( 1 - \rho) } {1-u} \frac{ f_{D}
(F^{-1} (u) ) } {f (F^{-1} (u) )
} \sqrt{\lambda} \frac{r_{D}} {r_{\bar{D}}} K_2(F_{\bar{D}}(F^{-1}(u)),
r_{\bar{D}})
\end{eqnarray*}
uniformly for $u \in[a, b]$, $r_{D} \in[c, 1]$
and $r_{\bar{D}} \in[d, 1 ]$ where $K_1$ and $K_2$ are
independent Kiefer processes.
%%%%%%%
\end{theorem}

%%%%%%%
The proof of Theorem \ref{PPVuseqproc} is complicated by the fact
that $\hat{S}_{D, r_{D}}(x)$\break and~$\hat{F}_{r_{D}, r_{\bar{D}}}^{-1}(t)$
are correlated because\vspace*{1pt} $\hat{F}_{r_{D},r_{\bar{D}}}(x)$ is a linear
combination\break of~$\hat{F}_{D,r_{D}}(x)$ and $\hat{F}_{\bar{D},r_{\bar
{D}}}(x)$. In contrast, the sequential empirical ROC curve and the
sequential empirical PPV curve indexed by the FPF are functionals of
two independent sequential empirical estimators, $\hat{S}_{D,r_{D}}(x)$
and $\hat{S}_{\bar{D},r_{\bar{D}}}^{-1}(t)$, which makes it easier to
show that $R_{r_{D}, r_{\bar{D}}}(t)$ and $P_{r_{D}, r_{\bar{D}}}(t)$
converge to the sum of independent Kiefer processes. To account for the
correlation between~$\hat{S}_{D, r_{D}}(x)$ and $\hat{F}_{r_{D}, r_{\bar
{D}}}^{-1}(t)$, we follow the approach of \citet{PandS68}, who
prove a
similar result for two correlated, fixed-sample empirical processes.
The proof of Theorem \ref{PPVuseqproc} can be found in the \hyperref[app]{Appendix}.

Theorem \ref{PPVuseqproc} also establishes asymptotic theory for the
sequential empirical NPV curve because $\widehat{\mathrm{NPV}}_{r_{D},r_{\bar
{D}}}(t)$ is a function of $\widehat{\mathrm{PPV}}_{r_{D}, r_{\bar{D}}}(t)$.
Corollary~\ref{NPVuseqproc} establishes the convergence of
$N_{r_{D},r_{\bar{D}}}(t)$ to the sum of two independent Kiefer processes.
%
%co3.7 #&#
%
\begin{corollary}
%%%%%%%
\label{NPVuseqproc}
%%%%%%%
Assume \textup{(A1)--(A4)} hold and let $\frac{ f_{D} (F^{-1} (u)
) } {f (F^{-1} (u) ) }$ be bounded on $
[a,b]$. As $n_{D} \rightarrow\infty$ and $n_{\bar{D}}
\rightarrow\infty$
\begin{eqnarray*}
N_{r_{D}, r_{\bar{D}}}(u) &\rightarrow_{d}&  - \frac{\rho( 1 -
\rho) } {u} \frac{ f_{\bar{D}} (F^{-1} (u)
) } {f (F^{-1} (u) ) } K_1(F_{D}(F^{-1}(u)),
r_{D}) \\
&&{} + \frac{\rho( 1 - \rho) } {u} \frac{ f_{D} (F^{-1}
(u) ) } {f (F^{-1} (u) ) } \sqrt
{\lambda} \frac{r_{D}} {r_{\bar{D}}} K_2(F_{\bar{D}}(F^{-1}(u)), r_{\bar{D}})
\end{eqnarray*}
uniformly for $u \in[a, b]$, $r_{D} \in[c, 1]$
and $r_{\bar{D}} \in[d, 1 ]$ where $K_1$ and $K_2$ are
independent Kiefer processes.
%%%%%%%
\end{corollary}

%%%%%%%
Corollary \ref{NPVuseqproc} is immediate from Theorem \ref
{PPVuseqproc} by noting that
\[
N_{r_{D}, r_{\bar{D}}}(t) = \frac{1-u}{u} P_{r_{D}, r_{\bar{D}}}(t).
\]
As with the ROC curve and the PPV curve indexed by the FPF, Theorem~%
\ref{PPVuseqproc} and Corollary \ref{NPVuseqproc} allow us to develop
distribution theory for summaries of the PPV and NPV curve indexed by
u. Distribution theory for a vector of points on the PPV or NPV curve
is left for the \hyperref[app]{Appendix} but we choose to highlight the joint
distribution of the sequential empirical estimate of a single point on
the PPV or NPV curve. Corollary \ref{PPVNPVuseq} establishes that
the sequential empirical estimate of a point on the PPV or NPV curve is
asymptotically normal and has independent increments when divided by
its variance.
%
%co3.8 #&#
%
\begin{corollary}
%%%%%%%
\label{PPVNPVuseq}
%%%%%%%
Assume \textup{(A1)--(A4)} hold and let $\frac{ f_{D} (F^{-1} (u)
) } {f (F^{-1} (u) ) }$ be bounded on $
[a,b]$. For $u \in(0, 1 )$ and J stopping times:
\begin{longlist}[(A)]
\item[(A)] $(\widehat{\mathrm{PPV}}_{r_{D,1}, r_{\bar{D},1}}(u), \widehat
{\mathrm{PPV}}_{r_{D,2}, r_{\bar{D},2}}(u), \ldots, \widehat{\mathrm{PPV}}_{r_{D,J},
r_{\bar{D},J}}(u) )$, is\vspace*{1pt}~approxima\-tely~multivariate normal with
\[
\widehat{\mathrm{PPV}}_{r_{D,i}, r_{\bar{D},i}}(u) \sim N \bigl(\mathrm{PPV}(u), \sigma
^{2}_{\widehat{\mathrm{PPV}}_{r_{D,i}, r_{\bar{D},i}}(u)} \bigr),\qquad
i = 1, 2, \ldots, J,
\]
and
\begin{eqnarray*}
&&\operatorname{Cov} [\widehat{\mathrm{PPV}}_{r_{D,i}, r_{\bar{D},i}}(u),
\widehat{\mathrm{PPV}}_{r_{D,j}, r_{\bar{D},j}}(u) ] \\
&&\qquad = \operatorname{Var} [ \widehat{\mathrm{PPV}}_{r_{D,j}, r_{\bar{D},j}}(u) ] =
\sigma^{2}_{\widehat{\mathrm{PPV}}_{r_{D,j}, r_{\bar{D},j}}(u)},\qquad
r_{i} \leq r_{j},
\end{eqnarray*}
where
\begin{eqnarray*}
\sigma^{2}_{\widehat{\mathrm{PPV}}_{r_{D,j}, r_{\bar{D},j}}(u)} & = &
\biggl(\frac{ f_{\bar{D}} (F^{-1} (u) ) } {f
(F^{-1} (u) ) } (1 - \rho) \biggr)^2
\operatorname{PPV}(u) \biggl( \frac{ \rho} {1 - u} - \mathrm{PPV}(u) \biggr) \\
&&{}\times\frac{1}{n_{D} r_{D,j}}  \\
&&{} + \biggl( \frac{ f_{D} (F^{-1} (u) ) }
{f (F^{-1} (u) ) } \rho\biggr)^2 \bigl( 1 -
\mathrm{PPV}(u) \bigr) \biggl( \frac{u - \rho} {1 - u} + \mathrm{PPV}(u) \biggr) \\
&&\hspace*{11pt}{}\times \frac{1}{
n_{\bar{D}} r_{\bar{D},j} } .
\end{eqnarray*}
\item[(B)] $(\widehat{\mathrm{NPV}}_{r_{D,1}, r_{\bar{D},1}}(u), \widehat
{\mathrm{NPV}}_{r_{D,2}, r_{\bar{D},2}}(u), \ldots, \widehat{\mathrm{NPV}}_{r_{D,J},
r_{\bar{D},J}}(u) )$, is~approxima\-tely~multivariate normal with,
\[
\widehat{\mathrm{NPV}}_{r_{D,i}, r_{\bar{D},i}}(u) \sim N \bigl(\mathrm{NPV}(u), \sigma
^{2}_{\widehat{\mathrm{NPV}}_{r_{D,i}, r_{\bar{D},i}}(u)} \bigr),\qquad
i = 1, 2, \ldots, J,
\]
and
\begin{eqnarray*}
&&\operatorname{Cov} [\widehat{\mathrm{NPV}}_{r_{D,i}, r_{\bar{D},i}}(u),
\widehat{\mathrm{NPV}}_{r_{D,j}, r_{\bar{D},j}}(u) ] \\
&&\qquad = \operatorname{Var} [ \widehat{\mathrm{NPV}}_{r_{D,j}, r_{\bar{D},j}}(u) ] =
\sigma^{2}_{\widehat{\mathrm{NPV}}_{r_{D,j}, r_{\bar{D},j}}(u)},\qquad
r_{i} \leq r_{j},
\end{eqnarray*}
where
\begin{eqnarray*}
\sigma^{2}_{\widehat{\mathrm{NPV}}_{r_{D,j}, r_{\bar{D},j}}(u)} & = &
\biggl(\frac{ f_{\bar{D}} (F^{-1} (u) ) } {f
(F^{-1} (u) ) } (1 - \rho) \biggr)^2
\biggl(\mathrm{NPV}(u) + \frac{\rho- u} {u} \biggr) \bigl( 1 - \mathrm{NPV}(u) \bigr)\\
&&{}\times
\frac{1}{n_{D} r_{D,j} } \\
&&{} + \biggl( \frac{ f_{D} (F^{-1} (u) ) }
{f (F^{-1} (u) ) } \rho\biggr)^2 \operatorname{NPV}(u) \biggl(
\frac{1-\rho}{u} - \mathrm{NPV}(u) \biggr)\\
&&\hspace*{11pt}{}\times \frac{1}{ n_{\bar{D}} r_{\bar{D},j} } .
\end{eqnarray*}
\end{longlist}
\end{corollary}

%%%%%%%
It is immediate from Theorem \ref{PPVuseqproc} and Corollary \ref
{NPVuseqproc} that $\widehat{\mathrm{PPV}}_{r_{D}, r_{\bar{D}}}(u)$ and
$\widehat{\mathrm{NPV}}_{r_{D}, r_{\bar{D}}}(u)$ are asymptotically normal with
an independent increments covariance structure. By noting that
\[
F_{D}(F^{-1}(u)) = 1 - \frac{1 - u} {\rho} \operatorname{PPV}(u) = \frac{u} {\rho}
\bigl( 1 - \mathrm{NPV}(u) \bigr)
\]
and
\[
F_{\bar{D}}(F^{-1}(u)) = 1 - \frac{1 - u} {1 - \rho} \bigl( 1 - \mathrm{PPV}(u)
\bigr) = \frac{u} {1 - \rho} \operatorname{NPV}(u),
\]
we can write the asymptotic variances of $\widehat{\mathrm{PPV}}_{r_{D}, r_{\bar
{D}}}(u)$ and $\widehat{\mathrm{NPV}}_{r_{D}, r_{\bar{D}}}(u)$ as functions of
$\mathrm{PPV}(u)$ and $\mathrm{NPV}(u)$, respectively. This provides a better
understanding of the mean-variance relationship for the asymptotic
distributions of $\widehat{\mathrm{PPV}}_{r_{D}, r_{\bar{D}}}(u)$ and $\widehat
{\mathrm{NPV}}_{r_{D}, r_{\bar{D}}}(u)$ and, perhaps, provides a form of the
variance that is easier to work with in practical situations (i.e.,
study design, estimating the standard error, etc.).

An important component of Theorem \ref{PPVuseqproc} and
Corollary \ref{NPVuseqproc} is that not only do $P_{r_{D}, r_{\bar
{D}}}(u)$ and $N_{r_{D}, r_{\bar{D}}}(u)$ converge to the sum of
independent Kiefer processes, but they both converge to the same two
Kiefer processes. As a~result, we are able to derive the correlation
between a point on the PPV curve and a point on the NPV curve.
Corollary \ref{PPVNPVujoint} provides a bivariate normal
approximation for a point on the PPV and a point on the NPV
curve.\looseness=-1

%co3.9 #&#
%
\begin{corollary}
%%%%%%%
\label{PPVNPVujoint}
%%%%%%%
Assume \textup{(A1)--(A4)} hold and let $\frac{ f_{D} (F^{-1} (u)
) } {f (F^{-1} (u) ) }$ be bounded on $
[a,b]$. For $u_{1}, u_{2} \in(0,1)$, $( \widehat
{\mathrm{PPV}}_{r_{D,1}, r_{\bar{D},1}}(u_{1}),\widehat{\mathrm{NPV}}_{r_{D,2}, r_{\bar
{D},2}}(u_{2}))$, is approximately bivariate normally distributed with
\[
\widehat{\mathrm{PPV}}_{r_{D, 1}, r_{\bar{D},1}}(u_{1}) \sim N\bigl(\mathrm{PPV}(u),
\sigma^{2}_{\widehat{\mathrm{PPV}}_{r_{D,1}, r_{\bar{D},1}}(u_{1})} \bigr)
\]
and
\[
\widehat{\mathrm{NPV}}_{r_{D, 2}, r_{\bar{D},2}}(u_{2}) \sim N\bigl(\mathrm{NPV}(u),
\sigma^{2}_{\widehat{\mathrm{NPV}}_{r_{D, 2}, r_{\bar{D},2}}(u_{2})}
\bigr)\vadjust{\goodbreak}
\]
with
\begin{eqnarray*}
&& \operatorname{Cov} [ \widehat{\mathrm{PPV}}_{r_{D, 1},
r_{\bar{D},1}}(u_{1}) , \widehat
{\mathrm{NPV}}_{r_{D, 2}, r_{\bar{D},2}}(u_{2}) ] \\[-2pt]
&&\qquad = \frac{ ( 1- \rho)^2 u_{1} (1 - u_{2} ) } {
(1 - u_{1} ) u_{2} }\frac{ f_{\bar{D}} (F^{-1}
(u_{1} ) ) } {f (F^{-1} (u_{1} ) ) }
\frac{ f_{\bar{D}} (F^{-1} (u_{2} ) ) } {f
(F^{-1} (u_{2} ) ) } \\[-2pt]
&&\qquad\quad{}\times\frac{ (r_{D,1} \wedge
r_{D,2} ) ( 1 - \mathrm{NPV}(u_{1}) ) \operatorname{PPV}(u_{2}
) }{n_{D} r_{D,1} r_{D,2} }\\[-2pt]
&&\qquad\quad{} + \frac{ \rho^2 u_{1} (1 - u_{2} ) } { (1 - u_{1}
) u_{2} }\frac{ f_{D} (F^{-1} (u_{1} ) ) }
{f (F^{-1} (u_{1} ) ) } \frac{ f_{D} (F^{-1}
(u_{2} ) ) } {f (F^{-1} (u_{2} )
) } \\[-2pt]
&&\qquad\quad\hspace*{11pt}{}\times\frac{ (r_{\bar{D},2} \wedge r_{\bar{D},2} )
\operatorname{NPV} (u_{1}) ( 1- \mathrm{PPV}(u_{2} ) ) }{n_{D} r_{D,1} r_{D,2} },
\end{eqnarray*}
when $u_{1} \leq u_{2}$ and
\begin{eqnarray*}
&& \operatorname{Cov} [ \widehat{\mathrm{PPV}}_{r_{D, 1}, r_{\bar{D},1}}(u_{1}) , \widehat
{\mathrm{NPV}}_{r_{D, 2}, r_{\bar{D},2}}(u_{2}) ] \\[-2pt]
&&\qquad = ( 1- \rho)^2 \frac{ f_{\bar{D}} (F^{-1}
(u_{1} ) ) } {f (F^{-1} (u_{1} ) ) }
\frac{ f_{\bar{D}} (F^{-1} (u_{2} ) ) } {f
(F^{-1} (u_{2} ) ) } \\[-2pt]
&&\qquad\quad{}\times\frac{(r_{D,1} \wedge
r_{D,2} ) ( 1 - \mathrm{NPV}(u_{2}) ) \operatorname{PPV}(u_{1}
) }{n_{D} r_{D,1} r_{D,2} }\\[-2pt]
&&\qquad\quad{} + \rho^2 \frac{ f_{D} (F^{-1} (u_{1} ) ) } {f
(F^{-1} (u_{1} ) ) } \frac{ f_{D} (F^{-1}
(u_{2} ) ) } {f (F^{-1} (u_{2} )
) } \\[-2pt]
&&\qquad\quad\hspace*{11pt}{}\times\frac{ (r_{\bar{D},2} \wedge r_{\bar{D},2} )
\operatorname{NPV}
(u_{2}) ( 1- \mathrm{PPV}(u_{1} ) ) }{n_{D} r_{D,1}
r_{D,2} },
\end{eqnarray*}
when $u_{2} \leq u_{1}$, where $ \sigma^{2}_{\widehat{\mathrm{PPV}}_{r_{D,1},
r_{\bar{D},1}}(u_{1})}$ and $\sigma^{2}_{\widehat{\mathrm{NPV}}_{r_{D, 2},
r_{\bar{D},2}}(u_{2})}$ are defined as in Corollary \ref{PPVNPVuseq}.
\end{corollary}

%%%%%%%
The case of a point on the PPV curve and a point on the NPV curve is
presented for simplicity but Corollary \ref{PPVNPVujoint} can be
extended to an arbitrary vector of points on the PPV and NPV curves.
Corollary \ref{PPVNPVujoint} has obvious practical implications. It
is not uncommon to classify the bottom $u_{1} \times100$\% of the
population as ``low-risk,'' the top $(1 - u_{2} ) \times
100$\% of the population as ``high-risk'' and the remainder of the
population as ``moderate-risk.'' In this case, one would be interested
in the NPV of the low-risk group and the PPV of the high-risk group.
Corollary \ref{PPVNPVujoint} provides the joint convergence of these
two estimates.

Finally, we note that asymptotic results for the fixed-sample empirical
PPV and NPV curves indexed by the percentile value of the marker
distribution can be derived as a special case of the results in this
section. It is immediate from Theorem \ref{PPVuseqproc} and
Corollary \ref{NPVuseqproc} that the fixed-sample empirical\vadjust{\goodbreak} PPV and
NPV curves converge to the sum of independent Brownian bridges by
letting $r_{D}$ and $r_{\bar{D}}$ both equal 1. Furthermore,
Corollary \ref{PPVNPVuseq} provides a normal approximation for the
fixed-sample empirical estimate of a~point on the PPV or NPV curve for
the special case when $J = 1$.\vspace*{-2pt}

%s4 #&#
\section{Finite sample properties}
\label{simsec}

\begin{figure}

\includegraphics{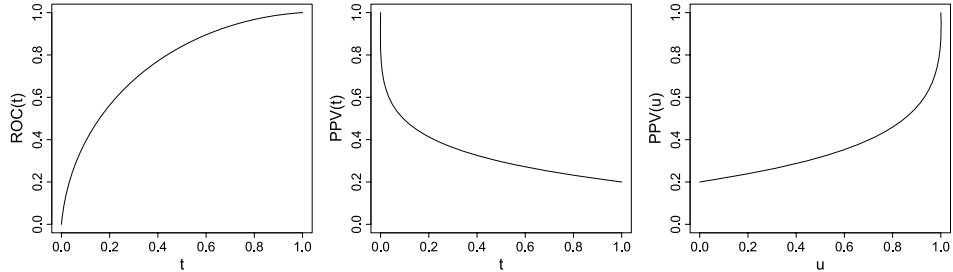}

\caption{True ROC and PPV curves for the scenario considered in
Section \protect\ref{simsec}.}
\label{simscenfig}
\vspace*{-2pt}
\end{figure}

A simulation study was completed to assess the finite sample properties
of the results in Theorems \ref{ROCseqproc}, \ref{PPVtseqproc} and
\ref{PPVuseqproc}. We simulated 10,000 studies with $n_{\bar{D}}$
controls and $n_{D}$ cases. Biomarker values for the controls were
drawn from a standard normal distribution and biomarker values for the
cases were drawn from a normal distribution with mean and standard
deviation equal to 1. A prevalence of 0.2 was used for estimation of the
PPV curve. Figure \ref{simscenfig} presents the true ROC and PPV
curves for this scenario. For each realization, we calculated
%
%f1 #&#
%
$R_{r_{D}, r_{\bar{D}}}(t)$, $P_{r_{D}, r_{\bar{D}}}
(t)$ and~$P_{r_{D}, r_{\bar{D}}}(u)$ and evaluated the
expected value, normality and covariance for various combinations of
$r_{D}$, $r_{\bar{D}}$ and $t$ or $u$. Normality was evaluated by
providing a summary of information found in a normal q-q plot. Instead
of providing the entire plot, we provide the (simulated) probability of
being less than the 5th, 25th, 50th, 75th and 95th percentile of a
normal distribution with variance derived using the results in
Theorems \ref{ROCseqproc}, \ref{PPVtseqproc} and \ref
{PPVuseqproc}. Similarly, the simulated covariance matrices were
compared to the theoretical covariance matrices derived using the
results in Theorems \ref{ROCseqproc}, \ref{PPVtseqproc} and~\ref
{PPVuseqproc}.

%
%t1 #&#
%
\begin{sidewaystable}
\tablewidth=\textheight
\tablewidth=\textwidth
\caption{Simulation results to evaluate the finite sample properties of
Theorem \protect\ref{ROCseqproc}. Presented are the expected value, simulated
probability of being less than 5th, 25th, 50th, 75th and 95th
percentile of the normal distribution, the simulated covariance matrix
and the theoretical covariance matrix for $R_{r_{D}, r_{\bar{D}}}
(t)$. 10,000 simulations were performed for each scenario}
%%%%%%%
\label{ROCsim}
\begin{tabular*}{\tablewidth}{@{\extracolsep{\fill}} l c c d{1.2} d{1.2}
d{1.2} c c d{1.3} d{1.3} d{1.3} d{1.3} c c c d{1.3} d{1.3}@{}}
\hline
& & \multicolumn{1}{c}{\textbf{5th}} & \multicolumn{1}{c}{\textbf{25th}}
& \multicolumn{1}{c}{\textbf{50th}} & \multicolumn{1}{c}{\textbf{75th}}
& \multicolumn{1}{c}{\textbf{95th}} && \multicolumn{4}{c}{\textbf{Observed}}
&&
\multicolumn{4}{c@{}}{\textbf{Theoretical}} \\
& \textbf{Mean} & \multicolumn{1}{c}{\textbf{\%tile}} & \multicolumn{1}{c}{\textbf{\%tile}}
& \multicolumn{1}{c}{\textbf{\%tile}} & \multicolumn{1}{c}{\textbf{\%tile}}
& \multicolumn{1}{c}{\textbf{\%tile}}
&& \multicolumn{4}{c}{\textbf{covariance matrix}}
&& \multicolumn{4}{c@{}}{\textbf{covariance matrix}} \\
\hline
\multicolumn{17}{@{}c@{}}{$n_{D} = 50$, $n_{\bar{D}} = 50$} \\
[4pt]
$R_{ 0.4 , 0.7 }( 0.4 )$ & 0.01 & 0.05 & 0.17 & 0.46 & 0.63 & 0.98 &&
0.1 & 0.117 & 0.079 & 0.103 & & 0.104 & 0.129 & 0.081 & 0.104 \\
$R_{ 1 , 1 }( 0.4 )$ & 0.02 & 0.07 & 0.2 & 0.44 & 0.74 & 0.97 & & & 0.318
& 0.104 & 0.262 & & & 0.322 & 0.104 & 0.26 \\
$R_{ 0.4 , 0.7 }( 0.2 )$ & 0.03 & 0.04 & 0.22 & 0.47 & 0.73 & 0.96 & &&
& 0.161 & 0.201 & & & & 0.171 & 0.225 \\
$R_{ 1 , 1 }( 0.2 )$ & 0.05 & 0.04 & 0.2 & 0.47 & 0.68 & 0.93 & & & &&
0.544 & & & & &0.563 \\
[4pt]
\multicolumn{17}{@{}c@{}}{$n_{D} = 100$, $n_{\bar{D}} = 100$} \\
[4pt]
$R_{ 0.4 , 0.7 }( 0.4 )$ & 0.01 & 0.05 & 0.21 & 0.41 & 0.78 & 0.97 &&
0.101 & 0.12 & 0.08 & 0.102 && 0.104 & 0.129 & 0.081 & 0.104 \\
$R_{ 1 , 1 }( 0.4 )$ & 0.02 & 0.05 & 0.24 & 0.48 & 0.76 & 0.96 & &&
0.318 & 0.104 & 0.26 & & & 0.322 & 0.104 & 0.26 \\
$R_{ 0.4 , 0.7 }( 0.2 )$ & 0.03 & 0.04 & 0.2 & 0.45 & 0.73 & 0.95 & &
&&
0.164 & 0.205 & & & &0.171 & 0.225 \\
$R_{ 1 , 1 }( 0.2 )$ & 0.05 & 0.04 & 0.23 & 0.47 & 0.73 & 0.95 & & & &&
0.55 & & & & &0.563 \\
[4pt]
\multicolumn{17}{@{}c@{}}{$n_{D} = 200$, $n_{\bar{D}} = 200$} \\
[4pt]
$R_{ 0.4 , 0.7 }( 0.4 )$ & 0.01 & 0.06 & 0.22 & 0.44 & 0.7 & 0.96 &&
0.104 & 0.121 & 0.081 & 0.102 && 0.104 & 0.129 & 0.081 & 0.104 \\
$R_{ 1 , 1 }( 0.4 )$ & 0.02 & 0.05 & 0.25 & 0.48 & 0.72 & 0.95 & &&
0.317 & 0.104 & 0.259 & & & 0.322 & 0.104 & 0.26 \\
$R_{ 0.4 , 0.7 }( 0.2 )$ & 0.03 & 0.04 & 0.25 & 0.5 & 0.7 & 0.94 & && &
0.168 & 0.212 & & & & 0.171 & 0.225 \\
$R_{ 1 , 1 }( 0.2 )$ & 0.05 & 0.05 & 0.23 & 0.46 & 0.72 & 0.95 & & && &
0.555 & & & & & 0.563 \\
\hline
\end{tabular*}
\end{sidewaystable}

Table \ref{ROCsim} presents\vspace*{1pt} simulation results for $R_{r_{D},r_{\bar
{D}}}(t)$. The expected value was close to 0 in all cases
with only a small amount of bias observed when $t = 0.2$. The
probability of being less than the theoretical 5th and 95th percentile
was close to the nominal value for all sample sizes, while the
probability of being less than the 25th, 50th and 75th percentile was
less than the nominal value with 50 cases and 50 controls but
approached the correct values as sample size increases. The observed
variance and covariance were less than expected with 50 cases and 50
controls but the observed covariance matrix approached the theoretical
covariance matrix in larger sample sizes. This phenomenon is likely due
to the sample space for $\mathrm{ROC}(t)$ being restricted to the unit interval.
$\widehat{\mathrm{ROC}}(t)$ is less likely to equal 0 or~1 as sample size
increases and the normal approximation will be more accurate. Similar
results were observed for $P_{r_{D}, r_{\bar{D}}}(t)$ and $P_{r_{D},
r_{\bar{D}}}(u)$ but were omitted for brevity.\vadjust{\goodbreak}

%s5 #&#
\section{Application}
%%%%%%%
\label{appsec}
%%%%%%%
The results of Section \ref{asympsec} provide fundamental theory that
allows existing group sequential methodology to be applied to summaries
of the ROC, PPV and NPV curves. In this section, we present an example
of how these results can be used to design group sequential diagnostic
biomarker studies. Our application is presented in the context of a
study to evaluate the diagnostic accuracy of des-gamma
carboxyprothrombin (DCP), a novel biomarker for the early detection of
hepatocellular carcinoma (HCC). A multi-center study was completed to
compare the diagnostic accuracy of DCP to that of alpha-fetoprotein
(AFP), the most widely used biomarker for the detection of HCC
[\citet
{MarrEtAl09}] but in our application we will only consider the design of
a study to compare DCP to historical levels of diagnostic accuracy for
AFP.

%%%%%%%
%%%%%%%
We consider a study to evaluate the predictive accuracy of DCP using
the following novel design that makes use of the joint asymptotic
theory for the PPV and NPV curve derived in Section
\ref{PPVNPVuseq}. Assume that the prevalence of HCC in the population of
interest is 0.2. In this case, one might call the bottom 60\% percent
of biomarker values ``negative,'' the top 10\% of the biomarker values
``positive'' and refer the remaining subjects for further evaluation.
Under this scenario, we would desire a high NPV for negative test
results, NPV(0.6), and a high PPV for positive test results, PPV(0.9).
The NPV(0.6) for AFP is 0.92 and the PPV(0.9) is 0.82. To determine if
DCP improves on the predictive accuracy of AFP, we would test the
hypothesis,
\[
H_{0}\mbox{:}\quad \mathrm{NPV}(0.6) \leq0.9
\quad\mbox{or}\quad\mathrm{PPV}(0.9) \leq0.8
\]
versus
\[
H_{a}\mbox{:}\quad \mathrm{NPV}(0.6) > 0.9
\quad\mbox{and}\quad\mathrm{PPV}(0.9) > 0.8
\]
using the test statistics, $Z_{\mathrm{NPV}(u_{1})}$ and $Z_{\mathrm{PPV}(u_{2})}$, where
$Z_{\mathrm{NPV}(u_{1})}$ is defined as
\[
Z_{\mathrm{NPV}(u_{1})} = \frac{ \widehat{\mathrm{NPV}}(0.6) - \mathrm{NPV}(0.6)_{0} } {\sigma_{\mathrm{NPV}(0.6)_{0}}},
\]
and $Z_{\mathrm{PPV}(u_{2})}$ is defined in an analogous fashion.

We consider a group sequential design using the error spending approach
proposed by \citet{HSandD90}. The overall null hypothesis will
only be
rejected if the null hypotheses for both $\mathrm{NPV}(0.6)$ and $\mathrm{PPV}(0.9)$ are
rejected. In the context of a group sequential study, this means that
the study will stop early to reject the null hypothesis if
$Z_{\mathrm{NPV}(u_{1})}$ and $Z_{\mathrm{PPV}(u_{2})}$ both cross the boundary for
rejecting the null hypothesis but the study will stop early for
futility if either $Z_{\mathrm{NPV}(u_{1})}$ or $Z_{\mathrm{PPV}(u_{2})}$ cross the
futility boundary. This implies that we do not need to adjust the
type-I error rate to account for multiple endpoints but we do need to
consider the joint probability of rejecting the null hypothesis when
determining the power.

The sample size for our study is chosen to achieve 90\% power under the
alternative hypothesis $\mathrm{NPV}(0.6) = 0.95$ and $\mathrm{PPV}(0.9) = 0.90$. A
closed-form formula for determining the required sample size is not
available. Instead, the sample size for a fixed sample design is
derived by numerically solving
\[
P\bigl(Z_{\mathrm{NPV}(u_{1})} > Z_{1-\alpha/2}, \mathrm{PPV}(u_{2}) > Z_{1-\alpha/2}
\vert \mathrm{NPV}(u_{1}) = 0.95, \mathrm{PPV}(u_{2}) = 0.90 \bigr)
\]
for $n_{D}$, where the joint distribution of $Z_{\mathrm{NPV}(u_{1})}$ and
$Z_{\mathrm{PPV}(u_{2})}$ is derived by applying the delta method to the joint
asymptotic normal distribution of~$\widehat{\mathrm{NPV}}_{r_{D}, r_{\bar
{D}}}(u_{1})$ and $\widehat{\mathrm{PPV}}_{r_{D}, r_{\bar{D}}}
(u_{2})$ found\vspace*{1pt} in Corollary \ref{PPVNPVujoint}. Assuming a
one-to-one ratio of cases to controls, 702 cases are required to
achieve 90\% power under the alternative hypothesis. This sample size
must be multiplied by an inflation factor to determine the maximum
sample size for a group sequential design (i.e., the sample size if the
study does not stop at the interim analyses) in order for the group
sequential design to maintain the same type-I error rate and power as
the fixed-sample design [\citet{JandT00}]. Using the gsDesign
package in
R, we find that the maximum sample size for group sequential studies
with two, three and four stopping times are 724, 737 and 745 cases,
respectively. However, as illustrated in the simulation which follows,
the actual sample sizes required in group sequential studies are
generally smaller than these maximum values.

%
%t2 #&#
%
\begin{table}
\tabcolsep=0pt
%%%%%%%
\caption{Simulation results to evaluate the operating characteristics
of a study to evaluate the predictive accuracy of DCP using a
fixed-sample design and group sequential designs with two, three or
four stopping times. Presented are the probability of rejecting the
null hypothesis and expected sample size under the null and alternative
hypotheses. 10,000 simulations were performed for each scenario}
%%%%%%%
\label{dcpppvnpvsim}
\begin{tabular*}{\tablewidth}{@{\extracolsep{\fill}}l c d{3.1} d{1.3} d{3.1} c d{3.1} c d{3.1}@{}}
\hline
& \multicolumn{2}{c}{$\bolds{\mathrm{NPV}(0.6) = 0.90}$}
& \multicolumn{2}{c}{$\bolds{\mathrm{NPV}(0.6) = 0.95}$}
& \multicolumn{2}{c}{$\bolds{\mathrm{NPV}(0.6) = 0.90}$}
& \multicolumn{2}{c@{}}{$\bolds{\mathrm{NPV}(0.6) = 0.95}$} \\
& \multicolumn{2}{c}{$\bolds{\mathrm{PPV}(0.9) = 0.80}$}
& \multicolumn{2}{c}{$\bolds{\mathrm{PPV}(0.9) = 0.80}$}
& \multicolumn{2}{c}{$\bolds{\mathrm{PPV}(0.9) =0.90}$}
& \multicolumn{2}{c@{}}{$\bolds{\mathrm{PPV}(0.9) = 0.90}$} \\[-4pt]
& \multicolumn{2}{c}{\hrulefill} & \multicolumn{2}{c}{\hrulefill}
& \multicolumn{2}{c}{\hrulefill} & \multicolumn{2}{c@{}}{\hrulefill}\\
\multirow{2}{35pt}[11pt]{\textbf{Stopping times}}
& \multicolumn{1}{c}{$\bolds{P(\mathrm{reject})}$}
& \multicolumn{1}{c}{$\bolds{E(n_{D})}$}
& \multicolumn{1}{c}{$\bolds{P(\mathrm{reject})}$}
& \multicolumn{1}{c}{$\bolds{E(n_{D})}$}
& \multicolumn{1}{c}{$\bolds{P(\mathrm{reject})}$}
& \multicolumn{1}{c}{$\bolds{E(n_{D})}$}
& \multicolumn{1}{c}{$\bolds{P(\mathrm{reject} )}$}
& \multicolumn{1}{c@{}}{$\bolds{E(n_{D})}$} \\
\hline
$J = 1$ & 0.003 & 702 & 0.03 & 702 & 0.026 & 702 & 0.917 & 702 \\
$J = 2$ & 0.004 & 432 & 0.026 & 492.4 & 0.024 & 489.5 & 0.924 & 624.5 \\
$J = 3$ & 0.004 & 367.4 & 0.022 & 431.3 & 0.023 & 433 & 0.917 & 580.1 \\
$J = 4$ & 0.002 & 340 & 0.023 & 410.7 & 0.024 & 417.2 & 0.911 & 571.1 \\
\hline
\end{tabular*}
\vspace*{-3pt}
\end{table}

Table \ref{dcpppvnpvsim} presents simulation results using a
fixed-sample design and group sequential designs with two, three and
four stopping times. Biomarker values for the controls were simulated
from a standard normal distribution and biomarker values for the cases
were simulated from a normal distribution with mean and variance chosen
to achieve the desired value of $\mathrm{NPV}(0.6)$ and $\mathrm{PPV}(0.9)$. The
advantages of group sequential designs are clear. The group sequential
designs have similar type-I error rate and power to the fixed-sample
design but with substantially smaller expected sample sizes in all
scenarios.\vadjust{\goodbreak}

%s6 #&#
\section{Discussion}
\label{discsec}
In this paper, we derived asymptotic properties of the sequential
empirical ROC, PPV and NPV curves. We first extended the work of
\citet
{HandT96} to the sequential empirical ROC curve and used these results
to develop distribution theory for summaries of the sequential
empirical ROC curve. Next, we considered asymptotic theory for the
sequential empirical PPV curve indexed by the FPF and percentile value
in the entire population. These results were used to develop
distribution theory for summaries of the sequential empirical PPV
curve. Asymptotic theory for the fixed-sample PPV curve, which was
previously unavailable, was developed as a special case.

This work was motivated by the desire to design group sequential
diagnostic biomarker studies. In Section \ref{appsec}, we illustrated
how our results can be used to design group sequential diagnostic
biomarker studies. Our simulation results clearly illustrate the
advantages of group sequential designs. In both cases, the group
sequential designs have similar type-I error rate and power than the
fixed-sample designs but with substantially smaller expected sample
size.

An advantage to our approach is that we are able to investigate the
joint behavior of multiple points on the ROC and PPV curve. The primary
endpoint of a diagnostic biomarker study may be a single point on the
ROC or PPV curve but other points on the ROC or PPV curve may also be
of interest. The results of Theorems \ref{ROCseqproc}, \ref
{PPVtseqproc} and \ref{PPVuseqproc} allow us to apply existing
group sequential methodology for analyzing multiple endpoints to
scenarios where multiple points on the ROC or PPV curve are of interest
in a group sequential diagnostic biomarker study [\citet{LandH01}].

We considered estimation of the sequential empirical ROC and PPV curve
under case-control sampling. The asymptotic properties of the
sequential empirical ROC and PPV curve under other sampling schemes are
also of interest. We are currently working on extending the results of
this paper to estimation of the sequential empirical ROC and PPV curve
under cohort and nested case-control sampling.

The theory developed in this paper applies to sequential testing of the
diagnostic accuracy of a continuous test. In many cases, diagnostic
tests take the form of multi-level ordinal data (cancer staging, for
example). Methods exist extending the ROC curve to ordinal data
[\citet
{DandA69}] but further work is needed to verify that group sequential
methods can be applied in these settings.

Response adaptive clinical trials have been proposed as a means to
provide greater flexibility when designing therapeutic clinical trials.
Response adaptive clinical trials adjust the design characteristics of
the study (sample size, percent randomized to each group, etc.) in
response to outcomes for subjects enrolled earlier in the study.
Recently, \citet{ZandH10} showed that a class of test statistics
from a
response adaptive clinical trial converges to Brownian Motion when
considered sequentially (similar to what we have shown for the emprical
ROC, PPV and NPV curves), which allows existing group sequential
methodology to be applied to response adaptive clinical trials.
Future
work will be needed to consider how response adaptive designs can be
applied in the setting of group sequential diagnostic biomarker
studies.\looseness=-1

\begin{appendix}\label{app}
%%%%%%%
% add "Appendix" to the section heading
%%%%%%%

%s7 #&#
\section*{\texorpdfstring{Appendix: Supplementary results for Section \lowercase{\protect\ref{asympsec}}}
{Appendix: Supplementary results for Section 3}}
%s7.1 #&#
\subsection{\texorpdfstring{Supplementary results for Section \protect\ref{ROCsec}}
{Supplementary results for Section 3.1}}

\mbox{}

\begin{pf*}{Proof of Theorem \ref{ROCseqproc}}
%%%%%%%
First, note that
\begin{eqnarray*}
&& n_{D}^{-1/2} [n_{D} r_{D}] \bigl(\widehat{\mathrm{ROC}}_{r_{D},
r_{\bar{D}}}(t) - \mathrm{ROC}(t)\bigr)\\
&&\qquad = n_{D}^{-1/2} [n_{D}
r_{D}] \bigl( \hat{S}_{D,r_{D}}(\hat{S}_{\bar
{D},r_{\bar{D}}}^{-1}(t)) - S_{D}(S_{\bar{D}}^{-1}(t)) \bigr) \\
&&\qquad = n_{D}^{-1/2} [n_{D} r_{D}] \bigl( \hat{S}_{D,r_{D}}(\hat{S}_{\bar
{D},r_{\bar{D}}}^{-1}(t)) - S_{D}(\hat{S}_{\bar{D},r_{\bar
{D}}}^{-1}(t)) \bigr) \\
&&\qquad\quad{} + n_{D}^{-1/2} [n_{D} r_{D}] \bigl( S_{D}(\hat{S}_{\bar{D},r_{\bar
{D}}}^{-1}(t)) - S_{D}(S_{\bar{D}}^{-1}(t)) \bigr).
\end{eqnarray*}
The first term converges to a Kiefer process. We note that
\begin{eqnarray*}
&&
\sup_{c \leq r_{D} \leq1} \sup_{d \leq r_{\bar{D}} \leq1} \sup_{ a
\leq t \leq b} | F_{\bar{D}} ( \hat{F}_{\bar{D}, r_{\bar
{D}}}^{-1} (t ) ) - t |\\
&&\qquad = \frac{ n_{\bar{D}}
} { [ n_{\bar{D}} d ]} \sup_{c \leq r_{D} \leq1} \sup_{d
\leq r_{\bar{D}} \leq1} \sup_{ a \leq t \leq b} \frac{ [ n_{\bar
{D}} d ]} { n_{\bar{D}} } | F_{\bar{D}} ( \hat{F}_{\bar
{D}, r_{\bar{D}}}^{-1} (t ) ) - t | \\
&&\qquad\leq \frac{ n_{\bar{D}} } { [ n_{\bar{D}} d ]} \sup_{c
\leq r_{D} \leq1} \sup_{d \leq r_{\bar{D}} \leq1} \sup_{ a \leq t
\leq b} \frac{ [ n_{\bar{D}} r_{\bar{D}} ]} { n_{\bar{D}} }
| F_{\bar{D}} ( \hat{F}_{\bar{D}, r_{\bar{D}}}^{-1} (t
) ) - t | .
\end{eqnarray*}
Therefore,
%
%e7.1 #&#
%
\begin{equation}
\label{gcroc}
{\sup_{c \leq r_{D} \leq1} \sup_{d \leq r_{\bar{D}} \leq1} \sup_{ a
\leq t \leq b} }| F_{\bar{D}} ( \hat{F}_{\bar{D}, r_{\bar
{D}}}^{-1} (t ) ) - t | \rightarrow_{\mathrm{a.s.}} 0
\end{equation}
by the Glivenko--Cantelli theorems [Theorems 1.51 and 1.52 in \citet
{CandS98}] and
because $\frac{ n_{\bar{D}} } { [ n_{\bar{D}} d ]}
\rightarrow\frac{1 } {d}$. Furthermore,\vspace*{1pt} $F_{\bar{D}}^{-1}(t)$ will be
continuous by \textup{(A1)--(A3)} and will be uniformly continuous on $[a,b]$. Therefore,
%
%e7.2 #&#
%
\begin{equation}
\label{rocquantconv}
{\sup_{c \leq r_{D} \leq1} \sup_{d \leq r_{\bar{D}} \leq1} \sup_{ a
\leq t \leq b}} | \hat{F}_{\bar{D}, r_{\bar{D}}}^{-1} (t
) - F_{\bar{D}}^{-1}(t) | \rightarrow_{\mathrm{a.s.}} 0.
\end{equation}
We note that due to the continuity of $F_{\bar{D}}(x)$, $S_{\bar
{D}}^{-1} (t ) = F_{\bar{D}}^{-1} (1 - t )$ and
therefore (\ref{rocquantconv}) also applies to $S_{\bar{D}}^{-1}
(t )$. From Corollary 1.A in \citet{CandS98}, (\ref
{rocquantconv}) and the uniform continuity of the Kiefer process, we have
%
%e7.3 #&#
%
\begin{equation}
\label{rocpart1}\qquad
n_{D}^{-1/2} [n_{D} r_{D}] \bigl( \hat{S}_{D,r_{D}}(\hat{S}_{\bar
{D},r_{\bar{D}}}^{-1}(t)) - S_{D}(\hat{S}_{\bar{D},r_{\bar
{D}}}^{-1}(t)) \bigr) \rightarrow_{d} K_1(\mathrm{ROC}(t),
r_{D}).\vadjust{\goodbreak}
\end{equation}
The second term can be rewritten as
\begin{eqnarray*}
&& n_{D}^{-1/2} [n_{D} r_{D}] \bigl( S_{D}(\hat{S}_{\bar{D},r_{\bar
{D}}}^{-1}(t)) - S_{D}(S_{\bar{D}}^{-1}(t)) \bigr) \\
%= & n_{D}^{-1/2} [n_{D} r_{D}] ( S_{D} ( S_{\bar{D}}^{-1}
%( S_{\bar{D}} ( \hat{S}_{\bar{D},r_{\bar{D}}}^{-1}(t)
%) ) ) - S_{D} ( S_{\bar{D}}^{-1}(t) )
%) \\
%= & \frac{n_{D}^{-1/2} [n_{D} r_{D}] } { n_{\bar{D}}^{-1/2} [n_{
%( S_{\bar{D}} ( \hat{S}_{\bar{D},r_{\bar{D}}}^{-1}(t)
%) ) ) - S_{D} ( S_{\bar{D}}^{-1}(t) )
%) } { S_{\bar{D}} ( \hat{S}_{\bar{D},r_{\bar{D}}}^{-1}(t)
%) - t } n_{\bar{D}}^{-1/2} [n_{\bar{D}} r_{\bar{D}}] ( S_{
%) \\
&&\qquad= \frac{n_{D}^{-1/2} [n_{D} r_{D}] } { n_{\bar{D}}^{-1/2} [n_{\bar
{D}} r_{\bar{D}}] } \frac{ ( S_{D} ( S_{\bar{D}}^{-1} (
S_{\bar{D}} ( \hat{S}_{\bar{D},r_{\bar{D}}}^{-1}(t) )
) ) - S_{D} ( S_{\bar{D}}^{-1}(t) ) ) } { S_{\bar
{D}} ( \hat{S}_{\bar{D},r_{\bar{D}}}^{-1}(t) ) - t } \\
&&\qquad\quad{}\times n_{\bar
{D}}^{-1/2} [n_{\bar{D}} r_{\bar{D}}] \bigl( S_{\bar{D}} ( \hat
{S}_{\bar{D},r_{\bar{D}}}^{-1}(t) ) - \hat{S}_{\bar{D},r_{\bar
{D}}} ( \hat{S}_{\bar{D},r_{\bar{D}}}^{-1}(t) ) \bigr) \\
&&\qquad\quad{} + \frac{n_{D}^{-1/2} [n_{D} r_{D}] } { n_{\bar{D}}^{-1/2} [n_{\bar
{D}} r_{\bar{D}}] } \frac{ ( S_{D} ( S_{\bar{D}}^{-1} (
S_{\bar{D}} ( \hat{S}_{\bar{D},r_{\bar{D}}}^{-1}(t) )
) ) - S_{D} ( S_{\bar{D}}^{-1}(t) ) ) } { S_{\bar
{D}} ( \hat{S}_{\bar{D},r_{\bar{D}}}^{-1}(t) ) - t }\\
&&\hspace*{11pt}\qquad\quad{}\times n_{\bar
{D}}^{-1/2} [n_{\bar{D}} r_{\bar{D}}] \bigl( \hat{S}_{\bar{D},r_{\bar
{D}}} ( \hat{S}_{\bar{D},r_{\bar{D}}}^{-1}(t) ) - t \bigr).
\end{eqnarray*}
By the mean value theorem, there exists a $S_{\bar{D}} ( \tilde
{S}_{\bar{D},r_{\bar{D}}}^{-1}(t) )$ between $S_{\bar{D}} (
\hat{S}_{\bar{D},r_{\bar{D}}}^{-1}(t) )$ and $t$ such that
\[
\frac{ S_{D} ( S_{\bar{D}}^{-1} ( S_{\bar{D}} ( \hat
{S}_{\bar{D},r_{\bar{D}}}^{-1}(t) ) ) ) - S_{D}
( S_{\bar{D}}^{-1}(t) ) } { S_{\bar{D}} ( \hat{S}_{\bar
{D},r_{\bar{D}}}^{-1}(t) ) - t } = \frac{f_{D} ( S_{\bar
{D}}^{-1} ( S_{\bar{D}} ( \tilde{S}_{\bar{D},r_{\bar
{D}}}^{-1}(t) ) ) ) } {f_{\bar{D}} ( S_{\bar
{D}}^{-1} ( S_{\bar{D}} ( \tilde{S}_{\bar{D},r_{\bar
{D}}}^{-1}(t) ) ) ) }.
\]
From (\ref{gcroc}), we know that $S_{\bar{D}} ( \hat{S}_{\bar
{D},r_{\bar{D}}}^{-1}(t) ) \rightarrow_{\mathrm{a.s.}} t$, uniformly for
$t \in[a,b]$, $r_{D} \in[c,1]$ and $r_{\bar{D}} \in[d,1]$, and,
therefore, $S_{\bar{D}} ( \tilde{S}_{\bar{D},r_{\bar{D}}}^{-1}(t)
) \rightarrow_{\mathrm{a.s.}} t$, uniformly for $t \in[a,b]$, $r_{D} \in
[c,1]$ and $r_{\bar{D}} \in[d,1]$. This, along with the uniform
continuity of $\frac{f_{D} ( S_{\bar{D}}^{-1} ( t )
) } {f_{\bar{D}} ( S_{\bar{D}}^{-1} ( t ) )
}$, allows us to conclude that
\[
\sup_{c \leq r_{D} \leq1} \sup_{d \leq r_{\bar{D}} \leq1} \sup_{ a
\leq t \leq b} \biggl| \frac{f_{D} ( S_{\bar{D}}^{-1} (
S_{\bar{D}} ( \tilde{S}_{\bar{D},r_{\bar{D}}}^{-1}(t) )
) ) } {f_{\bar{D}} ( S_{\bar{D}}^{-1} ( S_{\bar
{D}} ( \tilde{S}_{\bar{D},r_{\bar{D}}}^{-1}(t) ) )
) } - \frac{f_{D} ( S_{\bar{D}}^{-1} ( t )
) } {f_{\bar{D}} ( S_{\bar{D}}^{-1} ( t ) ) }
\biggr| \rightarrow_{\mathrm{a.s.}} 0,
\]
which implies
%
%e7.4 #&#
%
\begin{eqnarray}
\label{rocder}
&&\sup_{c \leq r_{D} \leq1} \sup_{d \leq r_{\bar{D}} \leq1} \sup_{ a
\leq t \leq b} \biggl| \frac{ S_{D} ( S_{\bar{D}}^{-1} (
S_{\bar{D}} ( \hat{S}_{\bar{D},r_{\bar{D}}}^{-1}(t) )
) ) - S_{D} ( S_{\bar{D}}^{-1}(t) ) } { S_{\bar{D}}
( \hat{S}_{\bar{D},r_{\bar{D}}}^{-1}(t) ) - t } \nonumber\\
&&\hspace*{199pt}{}- \frac
{f_{D} ( S_{\bar{D}}^{-1} ( t ) ) } {f_{\bar{D}}
( S_{\bar{D}}^{-1} ( t ) ) } \biggr| \\
&&\qquad\rightarrow
_{\mathrm{a.s.}} 0.\nonumber
\end{eqnarray}
For all $r_{\bar{D}} \in[d, 1]$,
\[
{\sup_{ a \leq t \leq b}} | \hat{S}_{\bar{D},r_{\bar{D}}} ( \hat
{S}_{\bar{D},r_{\bar{D}}}^{-1}(t) ) - t | \leq_{\mathrm{a.s.}} \frac
{ 1 } { [n_{\bar{D}} r_{\bar{D}}] }.
\]
Therefore,
\[
\sup_{c \leq r_{D} \leq1} \sup_{d \leq r_{\bar{D}} \leq1} \sup_{ a
\leq t \leq b} n_{\bar{D}}^{-1/2} [n_{\bar{D}} r_{\bar{D}}] | \hat
{S}_{\bar{D},r_{\bar{D}}} ( \hat{S}_{\bar{D},r_{\bar{D}}}^{-1}(t)
) - t | \leq_{\mathrm{a.s.}} \frac{ 1 } { n_{\bar{D}}^{1/2} }\vspace*{-2pt}
\]
and
%
%e7.5 #&#
%
\begin{equation}
\label{rochathat}
\sup_{c \leq r_{D} \leq1} \sup_{d \leq r_{\bar{D}} \leq1} \sup_{ a
\leq t \leq b} n_{\bar{D}}^{-1/2} [n_{\bar{D}} r_{\bar{D}}] | \hat
{S}_{\bar{D},r_{\bar{D}}} ( \hat{S}_{\bar{D},r_{\bar{D}}}^{-1}(t)
) - t | \rightarrow_{\mathrm{a.s.}} 0.\vspace*{-2pt}
\end{equation}
From Corollary 1.A in \citet{CandS98}, (\ref{rocquantconv}) and the
uniform continuity of the Kiefer process, we have
%
%e7.6 #&#
%
\begin{equation}
\label{rocpart2kief}
n_{\bar{D}}^{-1/2} [n_{\bar{D}} r_{\bar{D}}] \bigl( S_{\bar{D}} (
\hat{S}_{\bar{D},r_{\bar{D}}}^{-1}(t) ) - \hat{S}_{\bar{D},r_{\bar
{D}}} ( \hat{S}_{\bar{D},r_{\bar{D}}}^{-1}(t) ) \bigr)
\rightarrow_{d} K_2(t, r_{\bar{D}}).\vspace*{-2pt}
\end{equation}
By (\ref{rocder}), (\ref{rochathat}), (\ref{rocpart2kief}) and
noting that $\frac{n_{D}^{-1/2} [n_{D} r_{D}] } { n_{\bar{D}}^{-1/2}
[n_{\bar{D}} r_{\bar{D}}] } \rightarrow\lambda^{1/2} \frac{r_{D}}
{r_{\bar{D}}}$, we conclude that
%
%e7.7 #&#
%
\begin{eqnarray}
\label{rocpart2}
&&
n_{D}^{-1/2} [n_{D} r_{D}] \bigl( S_{D}(\hat{S}_{\bar{D},r_{\bar
{D}}}^{-1}(t)) - S_{D}(S_{\bar{D}}^{-1}(t)) \bigr)\nonumber\\[-9pt]\\[-9pt]
&&\qquad \rightarrow_{d}
\lambda^{1/2} \frac{r_{D}}{r_{\bar{D}}} \biggl(\frac{f_{D}(S_{\bar
{D}}^{-1}(t))} {f_{\bar{D}}(S_{\bar{D}}^{-1}(t)) } \biggr) K_2(t,
r_{\bar{D}}).\nonumber\vspace*{-2pt}
\end{eqnarray}
Summing (\ref{rocpart1}) and (\ref{rocpart2}) gives the desired result.\vspace*{-2pt}
%%%%%%%
\end{pf*}
%

%s7.2 #&#
\subsection{\texorpdfstring{Supplementary results for Section \protect\ref{PPVtsec}}
{Supplementary results for Section 3.2}}\vspace*{-2pt}
%
%co7.1 #&#
%
\begin{corollary}
%%%%%%%
\label{PPVmultt}
%%%%%%%
Assume \textup{(A1)--(A4)} hold and let $\frac{f_{D}(S_{\bar{D}}^{-1}(t))} {f_{\bar
{D}}(S_{\bar{D}}^{-1}(t)) }$ be bounded on $[a,b]$. For $t_{1}, t_{2},
\ldots, t_{J} \in(0,1)$, $r_{D,1}, r_{D,2}, \ldots, r_{D,J}
\in(0,1]$ and\vspace*{1pt} $r_{\bar{D},1}, r_{\bar{D},2},\allowbreak \ldots, r_{\bar
{D},J} \in(0,1]$, a vector\vspace*{1pt} of arbitrary points on the
sequential empirical PPV curve, $(\widehat{\mathrm{PPV}}_{r_{D,1}, r_{\bar
{D},1}}(t_{1}), \widehat{\mathrm{PPV}}_{r_{D,2}, r_{\bar{D},2}}(t_{2}), \ldots,
\widehat{\mathrm{PPV}}_{r_{D,J}, r_{\bar{D},J}}(t_{J}) )$, is
approximately multivariate normal with
\begin{eqnarray*}
\widehat{\mathrm{PPV}}_{r_{D,j}, r_{\bar{D},j}}(t_{j}) &\sim& N \bigl(\mathrm{PPV}(t_{j}),
\sigma^{2}_{\widehat{\mathrm{PPV}}_{r_{D,j}, r_{\bar{D},j}}(t_{j})}
\bigr),\qquad
j = 1, 2, \ldots, J,
\\[-2pt]
\sigma^{2}_{\widehat{\mathrm{PPV}}_{r_{D,j}, r_{\bar{D},j}}(t_{j})} &=&
\biggl(\frac{t (1 - \rho) \rho} { (\mathrm{ROC}(t) \rho+ t
( 1 - \rho) )^2 } \biggr)^2
\sigma^{2}_{\widehat{\mathrm{ROC}}_{r_{D,j}, r_{\bar{D},j}}(t_{j})}\vspace*{-2pt}
\end{eqnarray*}
and
\begin{eqnarray*}
&& \operatorname{Cov} [\widehat{\mathrm{PPV}}_{r_{D,i}, r_{\bar{D},i}}(t_{i}), \widehat
{\mathrm{PPV}}_{r_{D,j}, r_{\bar{D},j}}(t_{j}) ] \\[-2pt]
&&\qquad= \biggl(\frac{t_{i} (1 - \rho) \rho} { (\mathrm{ROC}(t_{i})
\rho+ t_{i} ( 1 - \rho) )^2 } \biggr) \biggl(\frac
{t_{j} (1 - \rho) \rho} { (\mathrm{ROC}(t_{j}) \rho+ t_{j}
( 1 - \rho) )^2 } \biggr)\\[-2pt]
&&\qquad\quad\hspace*{0pt}{}\times \operatorname{Cov} [\widehat
{\mathrm{ROC}}_{r_{D,i}, r_{\bar{D},i}}(t_{i}), \widehat{\mathrm{ROC}}_{r_{D,j}, r_{\bar
{D},j}}(t_{j}) ] ,\vspace*{-2pt}
\end{eqnarray*}
where $\sigma^{2}_{\widehat{\mathrm{ROC}}_{r_{D,j}, r_{\bar{D},j}}(t_{j})}$ and
$\operatorname{Cov} [\widehat{\mathrm{ROC}}_{r_{D,i}, r_{\bar{D},i}}(t_{i}), \widehat
{\mathrm{ROC}}_{r_{D,j}, r_{\bar{D},j}}(t_{j}) ] $ are as defined in
Corollary \ref{ROCmultt}.\vspace*{-2pt}\vadjust{\goodbreak}
%%%%%%%
\end{corollary}

\begin{pf}
%%%%%%%
Immediate from Theorem \ref{PPVtseqproc}.
%%%%%%%
\end{pf}
%
%co7.2 #&#
%
\begin{corollary}
%%%%%%%
\label{PPVtseq}
%%%%%%%
Assume\vspace*{1pt} \textup{(A1)--(A4)} hold and let $\frac{f_{D}(S_{\bar{D}}^{-1}(t))} {f_{\bar
{D}}(S_{\bar{D}}^{-1}(t)) }$ be bounded on $[a,b]$. For $t \in(0,1)$
and J stopping\vspace*{1pt} times $(\widehat{\mathrm{PPV}}_{r_{D,1}, r_{\bar{D},1}}(t),
\widehat{\mathrm{PPV}}_{r_{D,2}, r_{\bar{D},2}}(t),\allowbreak \ldots, \widehat
{\mathrm{PPV}}_{r_{D,J}, r_{\bar{D},J}}(t))$, is approximately
multivariate normal with
\[
\widehat{\mathrm{PPV}}_{r_{D,i}, r_{\bar{D},i}}(t) \sim N \bigl( \mathrm{PPV}(t), \sigma
^{2}_{\widehat{\mathrm{PPV}}_{r_{D,i}, r_{\bar{D},i}}(t)} \bigr),\qquad
i = 1, 2, \ldots, J,
\]
and
\begin{eqnarray*}
&&\operatorname{Cov} [ \widehat{\mathrm{PPV}}_{r_{D,i}, r_{\bar{D},i}}(t),
\widehat{\mathrm{PPV}}_{r_{D,j}, r_{\bar{D},j}}(t) ] \\
&&\qquad= \operatorname{Var} [ \widehat{\mathrm{PPV}}_{r_{D,j}, r_{\bar{D},j}}(t) ] = \sigma
^{2}_{\widehat{\mathrm{PPV}}_{r_{D,j}, r_{\bar{D},j}}(t)},\qquad
r_{i}\leq r_{j},
\end{eqnarray*}
where $\sigma^{2}_{\widehat{\mathrm{PPV}}_{r_{D,j}, r_{\bar{D},j}}(t)}$ is
defined as in Corollary \ref{PPVmultt}.
%%%%%%%
\end{corollary}
\begin{pf}
%%%%%%%
Immediate from Corollary \ref{PPVmultt}.
%%%%%%%
\end{pf}
%

%s7.3 #&#
\subsection{\texorpdfstring{Supplementary results for Section \protect\ref{PPVusec}}
{Supplementary results for Section 3.3}}

\mbox{}

\begin{pf*}{Proof of Theorem \ref{PPVuseqproc}}
%%%%%%
The proof of Theorem \ref{PPVuseqproc} follows the proofs found in
\citet{PandS68}. First, note that
\begin{eqnarray*}
&&n_{D}^{-1/2} [ n_{D} r_{D} ] \bigl( \hat{S}_{D,r_{D}} (
\hat{F}_{r_{D}, r_{\bar{D}}}^{-1} (u) )
- S_{D} (F^{-1} (u) ) \bigr) \\
%= & n_{D}^{-1/2} [ n_{D} r_{D} ] ( F_{D} (F^{-1}
%(u) ) - \hat{F}_{D,r_{D}} ( \hat{F}_{r_{D}, r_{
&&\qquad= n_{D}^{-1/2} [ n_{D} r_{D} ] \bigl( F_{D} (F^{-1}
(u) ) - F_{D} ( \hat{F}_{r_{D}, r_{\bar{D}}}^{-1}
(u) ) \bigr) \\
&&\qquad\quad{} + n_{D}^{-1/2} [ n_{D} r_{D} ] \bigl( F_{D} ( \hat
{F}_{r_{D}, r_{\bar{D}}}^{-1} (u) ) - \hat
{F}_{D,r_{D}} ( \hat{F}_{r_{D}, r_{\bar{D}}}^{-1} (u)
) \bigr).
\end{eqnarray*}
The first term can be rewritten as
\begin{eqnarray*}
&& n_{D}^{-1/2} [ n_{D} r_{D} ] \bigl( F_{D} (F^{-1}
(u) ) - F_{D} ( \hat{F}_{r_{D}, r_{\bar{D}}}^{-1}
(u) ) \bigr) \\
%= & n_{D}^{-1/2} [ n_{D} r_{D} ] \frac{ F_{D} (F^{-1}
%( u ) ) - F_{D} ( F^{-1} ( F (
%) } { u - F ( \hat{F}_{r_{D}, r_{\bar{D}}}^{-1} (u
%) ) } ( u - F ( \hat{F}_{r_{D}, r_{
%= & n_{D}^{-1/2} [ n_{D} r_{D} ] \frac{ F_{D} (
%F^{-1} ( F ( \hat{F}_{r_{D}, r_{\bar{D}}}^{-1} (u
%) ) ) ) - F_{D} (F^{-1} ( u )
%) } { F ( \hat{F}_{r_{D}, r_{\bar{D}}}^{-1} (u)
%) - u } ( u - \hat{F}_{r_{D}, r_{\bar{D}}} (
%& + n_{D}^{-1/2} [ n_{D} r_{D} ] \frac{ F_{D} (
%F^{-1} ( F ( \hat{F}_{r_{D}, r_{\bar{D}}}^{-1} (u
%) ) ) ) - F_{D} (F^{-1} ( u )
%) } { F ( \hat{F}_{r_{D}, r_{\bar{D}}}^{-1} (u)
%) - u } ( \hat{F}_{r_{D}, r_{\bar{D}}} (
&&\qquad= \frac{ F_{D} ( F^{-1} ( F ( \hat{F}_{r_{D}, r_{\bar
{D}}}^{-1} (u) ) ) ) - F_{D} (F^{-1}
( u ) ) } { F ( \hat{F}_{r_{D}, r_{\bar{D}}}^{-1}
(u) ) - u }\\
&&\qquad\quad{}\times n_{D}^{-1/2} [ n_{D} r_{D} ]
\bigl( u - \hat{F}_{r_{D}, r_{\bar{D}}} ( \hat{F}_{r_{D}, r_{\bar
{D}}}^{-1} (u) ) \bigr) \\
&&\qquad\quad{} + \frac{ F_{D} ( F^{-1} ( F ( \hat{F}_{r_{D}, r_{\bar
{D}}}^{-1} (u) ) ) ) - F_{D} (F^{-1}
( u ) ) } { F ( \hat{F}_{r_{D}, r_{\bar{D}}}^{-1}
(u) ) - u } \\
&&\hspace*{11pt}\qquad\quad{}\times\rho n_{D}^{-1/2} [ n_{D} r_{D}
] \bigl( \hat{F}_{D, r_{D}} ( \hat{F}_{r_{D}, r_{\bar{D}}}^{-1}
(u) ) - F_{D} ( \hat{F}_{r_{D}, r_{\bar{D}}}^{-1}
(u) ) \bigr) \\
&&\qquad\quad{} + \frac{n_{D}^{-1/2} [ n_{D} r_{D} ]} { n_{\bar
{D}}^{-1/2} [ n_{\bar{D}} r_{\bar{D}} ] }\frac{ F_{D}
( F^{-1} ( F ( \hat{F}_{r_{D}, r_{\bar{D}}}^{-1} (u
) ) ) ) - F_{D} (F^{-1} ( u ) )
} { F ( \hat{F}_{r_{D}, r_{\bar{D}}}^{-1} (u) ) -
u } \\
&&\hspace*{11pt}\qquad\quad{}\times(1 - \rho) n_{\bar{D}}^{-1/2} [ n_{\bar{D}} r_{\bar
{D}} ] \bigl( \hat{F}_{\bar{D}, r_{\bar{D}}} ( \hat
{F}_{r_{D}, r_{\bar{D}}}^{-1} (u) ) - F_{\bar{D}}
( \hat{F}_{r_{D}, r_{\bar{D}}}^{-1} (u) ) \bigr).
\end{eqnarray*}
We begin by showing that $F ( \hat{F}_{r_{D}, r_{\bar{D}}}^{-1}
(u) )$ converges to $u$ uniformly,
\begin{eqnarray*}
&&{\sup_{c \leq r_{D} \leq1} \sup_{ d \leq r_{\bar{D}} \leq1 }  \sup_{
a \leq u \leq b } }| F ( \hat{F}_{r_{D}, r_{\bar{D}}}^{-1}
(u) ) - u | \\
&&\qquad\leq {\sup_{c \leq r_{D} \leq1} \sup_{ d \leq r_{\bar{D}} \leq1 }
\sup_{ a \leq u \leq b }} | F ( \hat{F}_{r_{D}, r_{\bar
{D}}}^{-1} (u) ) - \hat{F}_{r_{D}, r_{\bar{D}}} (
\hat{F}_{r_{D}, r_{\bar{D}}}^{-1} (u) ) | \\
&&\qquad\quad{} + {\sup_{c \leq r_{D} \leq1} \sup_{ d \leq r_{\bar{D}} \leq1 } \sup
_{ a \leq u \leq b }} | \hat{F}_{r_{D}, r_{\bar{D}}} ( \hat
{F}_{r_{D}, r_{\bar{D}}}^{-1} (u) ) - u |.
\end{eqnarray*}
We note that
\begin{eqnarray*}
&&{\sup_{c \leq r_{D} \leq1} \sup_{ d \leq r_{\bar{D}} \leq1 } \sup_{
a \leq u \leq b }} | F ( \hat{F}_{r_{D}, r_{\bar{D}}}^{-1}
(u) ) - \hat{F}_{r_{D}, r_{\bar{D}}} ( \hat
{F}_{r_{D}, r_{\bar{D}}}^{-1} (u) ) | \\
%& + \sup_{c \leq r_{D} \leq1} \sup_{ d \leq r_{\bar{D}} \leq1 }
%( \hat{F}_{r_{D}, r_{\bar{D}}}^{-1} (u) )
%| \\
%= & \frac{ n_{D} } { [n_{D}c] } \sup_{c \leq r_{D} \leq1} \sup_{ d
%{ n_{D} } | F_{D} ( \hat{F}_{r_{D}, r_{\bar{D}}}^{-1} (u
%) ) - \hat{F}_{D, r_{D}} ( \hat{F}_{r_{D}, r_{
%& + \frac{ n_{\bar{D}} } { [n_{\bar{D}}d] } \sup_{c \leq r_{D} \leq
%1} \sup_{ d \leq r_{\bar{D}} \leq1 } \sup_{ a \leq u \leq b } \frac{
%[n_{\bar{D}}d] } { n_{\bar{D}} } | F_{\bar{D}} (
%) ) | \\
&&\quad\leq \frac{ n_{D} } { [n_{D}c] } \sup_{c \leq r_{D} \leq1} \sup_{ d
\leq r_{\bar{D}} \leq1 } \sup_{ a \leq u \leq b } \frac{ [n_{D}r_{D}]
} { n_{D} } | F_{D} ( \hat{F}_{r_{D}, r_{\bar{D}}}^{-1}
(u) ) - \hat{F}_{D, r_{D}} ( \hat{F}_{r_{D}, r_{\bar
{D}}}^{-1} (u) ) | \\
&&\qquad + \frac{ n_{\bar{D}} } { [n_{\bar{D}}d] } \sup_{c \leq r_{D} \leq1}
\sup_{ d \leq r_{\bar{D}} \leq1 } \sup_{ a \leq u \leq b } \frac{
[n_{\bar{D}}r_{\bar{D}}] } { n_{\bar{D}} } | F_{\bar{D}} (
\hat{F}_{r_{D}, r_{\bar{D}}}^{-1} (u) )\\
&&\hspace*{207pt}{} - \hat{F}_{\bar
{D}, r_{\bar{D}}} ( \hat{F}_{r_{D}, r_{\bar{D}}}^{-1} (u
) ) | \\
&&\quad \rightarrow_{\mathrm{a.s.}} 0,
\end{eqnarray*}
by the Glivenko--Cantelli theorems [Theorems 1.51 and 1.52 in \citet{CandS98}],
along with the fact that $\frac{n_{D}}{[n_{D} c]} \rightarrow\frac
{1} {c}$ and $\frac{n_{\bar{D}}}{[n_{\bar{D}} d]} \rightarrow\frac
{1} {d}$. For all $r_{D}, r_{\bar{D}} \in(0, 1 ] \times
(0, 1 ]$,
\[
\sup_{ a \leq u \leq b } | u - \hat{F}_{r_{D}, r_{\bar{D}}} (
\hat{F}_{r_{D}, r_{\bar{D}}}^{-1} (u) ) | \leq
_{\mathrm{a.s.}} \biggl( \frac{\rho} {[ r_{D} n_{D} ] } \vee
\frac{ 1 - \rho} { [ n_{\bar{D}} r_{\bar{D}} ] } \biggr).
\]
Therefore,
\begin{eqnarray*}
&&\sup_{c \leq r_{D} \leq1} \sup_{ d \leq r_{\bar{D}} \leq1 } \sup_{ a
\leq u \leq b } | u - \hat{F}_{r_{D}, r_{\bar{D}}} ( \hat
{F}_{r_{D}, r_{\bar{D}}}^{-1} (u) ) | \\
&&\qquad\leq
_{\mathrm{a.s.}} \biggl( \frac{\rho} {[ n_{D} c ] } \vee\frac{
1 - \rho} { [ n_{\bar{D}} d ] } \biggr) \rightarrow0,
\end{eqnarray*}
which implies that
%
%e7.8 #&#
%
\begin{equation}
\label{FFinvtov}
\sup_{c < r_{D} \leq1} \sup_{d < r_{\bar{D}} \leq1 } \sup_{ a \leq u
\leq b } | F ( \hat{F}_{r_{D}, r_{\bar{D}}}^{-1} (u
) ) - u | \rightarrow_{\mathrm{a.s.}} 0.
\end{equation}
We note that (\ref{FFinvtov}) also implies that $F_{D} ( \hat
{F}_{r_{D}, r_{\bar{D}}}^{-1} (u) )$ and $F_{\bar{D}}
( \hat{F}_{r_{D}, r_{\bar{D}}}^{-1} (u) )$
converge uniformly\vadjust{\goodbreak} to $F_{D} ( F^{-1} (u) )$ and
$F_{\bar{D}} ( F^{-1} (u) )$, respectively, which
can be seen by noting that the difference between $F_{D} ( \hat
{F}_{r_{D}, r_{\bar{D}}}^{-1} (u) )$ and $F_{D} (
F^{-1} (u) )$ will always have the same sign as the
difference between $F_{\bar{D}} ( \hat{F}_{r_{D}, r_{\bar
{D}}}^{-1} (u) )$ and~$F_{\bar{D}} ( F^{-1}
(u) )$.

By the mean value theorem, there exists $F ( \tilde{F}_{r_{D},
r_{\bar{D}}}^{-1} (u) )$ between $u$\break and~$F ( \hat
{F}_{r_{D}, r_{\bar{D}}}^{-1} (u) )$, such that
\[
\frac{ F_{D} ( F^{-1} ( F ( \hat{F}_{r_{D}, r_{\bar
{D}}}^{-1} (u) ) ) ) - F_{D} (F^{-1}
( u ) ) } { F ( \hat{F}_{r_{D}, r_{\bar{D}}}^{-1}
(u) ) - u } = \frac{ f_{D} ( F^{-1} ( F
( \tilde{F}_{r_{D}, r_{\bar{D}}}^{-1} (u) ) )
) } { f ( F^{-1} ( F ( \tilde{F}_{r_{D}, r_{\bar
{D}}}^{-1} (u) ) ) ) }.
\]
The uniform continuity of $\frac{ f_{D} ( F^{-1} (u)
) } { f ( F^{-1} (u) ) }$, combined with
the fact that
\[
F ( \tilde{F}_{r_{D}, r_{\bar{D}}}^{-1}
(u) ) \rightarrow_{\mathrm{a.s.}} u
\]
uniformly, allows us to conclude
%
%e7.9 #&#
%
\begin{equation}
\label{der}
\sup_{c < r_{D} \leq1} \sup_{d < r_{\bar{D}} \leq1 }
\sup_{ a \leq u \leq b } \biggl| \frac{ f_{D} ( F^{-1} ( F
( \tilde{F}_{r_{D}, r_{\bar{D}}}^{-1} (u) )
) ) } { f ( F^{-1} ( F ( \tilde{F}_{r_{D}, r_{\bar
{D}}}^{-1} (u) ) ) ) } - \frac{ f_{D}
( F^{-1} (u) ) } { f ( F^{-1} (u)
) } \biggr| \rightarrow_{\mathrm{a.s.}} 0.\hspace*{-32pt}
\end{equation}
For all $r_{D}, r_{\bar{D}} \in(0, 1 ] \times(0, 1
]$,
\begin{eqnarray*}
&&\sup_{ a \leq u \leq b } n_{D}^{-1/2} [ n_{D} r_{D} ] |
u - \hat{F}_{r_{D}, r_{\bar{D}}} ( \hat{F}_{r_{D}, r_{\bar
{D}}}^{-1} (u) ) | \\
&&\qquad\leq_{\mathrm{a.s.}} \biggl( \frac
{\rho} { n_{D}^{-1/2} } \vee\frac{[ n_{D} r_{D} ] } {
[ n_{\bar{D}} r_{\bar{D}} ] } \frac{ 1 - \rho} {
n_{D}^{-1/2} } \biggr).
\end{eqnarray*}
Therefore, as $n_{D} \rightarrow\infty$ and $n_{\bar{D}} \rightarrow
\infty$,
\[
\sup_{ 0 < r_{D} \leq1} \sup_{ 0 < r_{\bar{D}} \leq1 } \sup_{ a \leq
u \leq b } n_{D}^{-1/2} [ n_{D} r_{D} ] | u - \hat
{F}_{r_{D}, r_{\bar{D}}} ( \hat{F}_{r_{D}, r_{\bar{D}}}^{-1}
(u) ) | \rightarrow_{\mathrm{a.s.}} 0.
\]
Combining this result with (\ref{der}) allows us to conclude that
\begin{eqnarray*}
&&\frac{ F_{D} ( F^{-1} ( F ( \hat{F}_{r_{D}, r_{\bar
{D}}}^{-1} (u) ) ) ) - F_{D} (F^{-1}
( u ) ) } { F ( \hat{F}_{r_{D}, r_{\bar{D}}}^{-1}
(u) ) - u } n_{D}^{-1/2} [ n_{D} r_{D} ]
\bigl( u - \hat{F}_{r_{D}, r_{\bar{D}}} ( \hat{F}_{r_{D}, r_{\bar
{D}}}^{-1} (u) ) \bigr)\\
&&\qquad \rightarrow_{\mathrm{a.s.}} 0.
\end{eqnarray*}
Corollary 1.A in \citet{CandS98}, (\ref{der}) and the uniform continuity
of the Kiefer process allow us to conclude
%
%e7.10 #&#
%
\begin{eqnarray} \label{secondterm}
&&\frac{ F_{D} ( F^{-1} ( F ( \hat{F}_{r_{D}, r_{\bar
{D}}}^{-1} (u) ) ) ) - F_{D} (F^{-1}
( u ) ) } { F ( \hat{F}_{r_{D}, r_{\bar{D}}}^{-1}
(u) ) - u }  \nonumber\\
&&\quad{}\times\rho n_{D}^{-1/2} [ n_{D} r_{D}
] \bigl( \hat{F}_{D, r_{D}} ( \hat{F}_{r_{D}, r_{\bar
{D}}}^{-1} (u) ) - F_{D} ( \hat{F}_{r_{D}, r_{\bar
{D}}}^{-1} (u) ) \bigr) \\
&&\qquad \rightarrow_{d} \frac{f_{D} (F^{-1} (u) ) } {
f (F^{-1} (u) ) } \rho K_{1} (F_{D} (
F^{-1} (u) ), r_{D} )\nonumber
\end{eqnarray}
and
%
%e7.11 #&#
%
\begin{eqnarray}\label{thirdterm}
&&\frac{n_{D}^{-1/2} [ n_{D} r_{D} ]} { n_{\bar{D}}^{-1/2}
[ n_{\bar{D}} r_{\bar{D}} ] }\frac{ F_{D} ( F^{-1}
( F ( \hat{F}_{r_{D}, r_{\bar{D}}}^{-1} (u)
) ) ) - F_{D} (F^{-1} ( u ) ) } { F
( \hat{F}_{r_{D}, r_{\bar{D}}}^{-1} (u) ) - u }\nonumber\\[-3pt]
&&\quad{}\times(1 - \rho) n_{\bar{D}}^{-1/2} [ n_{\bar{D}} r_{\bar
{D}} ] \bigl( \hat{F}_{\bar{D}, r_{\bar{D}}} ( \hat
{F}_{r_{D}, r_{\bar{D}}}^{-1} (u) ) - F_{\bar{D}}
( \hat{F}_{r_{D}, r_{\bar{D}}}^{-1} (u) ) \bigr)
\\[-3pt]
&&\qquad \rightarrow_{d} \sqrt{\lambda} \frac{r_{D}} {r_{\bar{D}}} \frac
{f_{D} (F^{-1} (u) ) } { f (F^{-1}
(u) ) } (1 - \rho) K_{2} (F_{\bar{D}}
( F^{-1} (u) ), r_{D} ).\nonumber\vspace*{-2pt}
\end{eqnarray}
The second term converges in distribution to a Kiefer process
%
%e7.12 #&#
%
\begin{eqnarray} \label{firstterm}
&&
n_{D}^{-1/2} [ n_{D} r_{D} ] \bigl( F_{D} ( \hat
{F}_{r_{D}, r_{\bar{D}}}^{-1} (u) ) - \hat
{F}_{D,r_{D}} ( \hat{F}_{r_{D}, r_{\bar{D}}}^{-1} (u)
) \bigr)\nonumber\\[-3pt]
&&\qquad=  - n_{D}^{-1/2} [ n_{D} r_{D} ] \bigl(
\hat{F}_{D,r_{D}} ( \hat{F}_{r_{D}, r_{\bar{D}}}^{-1} (u
) ) - F_{D} ( \hat{F}_{r_{D}, r_{\bar{D}}}^{-1} (u
) ) \bigr) \\[-3pt]
%= & - n_{D}^{-1/2} [ n_{D} r_{D} ] ( \hat{F}_{D,r_{D}}
%( F_{D}^{-1} ( F_{D} ( \hat{F}_{r_{D}, r_{
%F_{D}^{-1} ( F_{D} ( \hat{F}_{r_{D}, r_{\bar{D}}}^{-1}
%( u ) ) ) ) ) \nonumber\\
&&\qquad\rightarrow_{d} - K_{1} ( F_{D} ( F^{-1} (u)
), r_{D} )\nonumber\vspace*{-2pt}
\end{eqnarray}
by Corollary 1.A in \citet{CandS98}. Summing (\ref{secondterm}),
(\ref
{thirdterm}) and (\ref{firstterm}) gives the desired result.\vspace*{-3pt}
%%%%%%%
\end{pf*}
%
%co7.3 #&#
%
\begin{corollary}
%%%%%%%
\label{PPVmultu}
%%%%%%%
Assume \textup{(A1)--(A4)} hold and let $\frac{ f_{D} (F^{-1} (u) ) }
{f (F^{-1} (u) ) }$ be bounded on $ [a,b]$. For $u_{1}, u_{2}, \ldots,
u_{J} \in(0,1)$, $r_{D,1}, r_{D,2}, \ldots, r_{D,J} \in(0,1]$ and
$r_{\bar {D},1}, r_{\bar{D},2},\allowbreak \ldots, r_{\bar{D},J} \in(0,1]$, a
vector of arbitrary points on the sequential\vspace*{1pt} empirical PPV curve, $
(\widehat{\mathrm{PPV}}_{r_{D,1}, r_{\bar{D},1}}(u_{1}),
\widehat{\mathrm{PPV}}_{r_{D,2}, r_{\bar{D},2}}(u_{2}), \ldots,
\widehat{\mathrm{PPV}}_{r_{D,J}, r_{\bar {D},J}}(u_{J}) )$, is
approximately multivariate normal with
\[
\widehat{\mathrm{PPV}}_{r_{D,j}, r_{\bar{D},j}}(u_{j}) \sim N \bigl(\mathrm{PPV}(u_{j}),
\sigma^{2}_{\widehat{\mathrm{PPV}}_{r_{D,j}, r_{\bar{D},j}}(u_{j})} \bigr),\qquad
j = 1, 2, \ldots, J,\vspace*{-2pt}
\]
with
\begin{eqnarray*}
&& \operatorname{Cov} [ \widehat{\mathrm{PPV}}_{r_{D, 1}, r_{\bar{D},1}}(u_{1}), \widehat
{\mathrm{PPV}}_{r_{D, 2}, r_{\bar{D},2}}(u_{2}) ] \\[-3pt]
&&\qquad = \frac{ ( 1- \rho)^2 u_{1} } { (1 - u_{1} )
}\frac{ f_{\bar{D}} (F^{-1} (u_{1} ) ) } {f
(F^{-1} (u_{1} ) ) } \frac{ f_{\bar{D}} (F^{-1}
(u_{2} ) ) } {f (F^{-1} (u_{2} )
) } \\[-3pt]
&&\qquad\quad{}\times\frac{ (r_{D,1} \wedge r_{D,2} ) ( 1 - \mathrm{NPV}
(u_{1}) ) \operatorname{PPV}(u_{2} ) }{n_{D} r_{D,1} r_{D,2} }\\[-3pt]
&&\qquad\quad{} + \frac{ \rho^2 u_{1} } { (1 - u_{1} ) }\frac{ f_{D}
(F^{-1} (u_{1} ) ) } {f (F^{-1} (u_{1}
) ) } \frac{ f_{D} (F^{-1} (u_{2} ) ) } {f
(F^{-1} (u_{2} ) ) } \\[-3pt]
&&\hspace*{11pt}\qquad\quad{}\times\frac{ (r_{\bar{D},2}
\wedge r_{\bar{D},2} ) \operatorname{NPV}(u_{1}) ( 1- \mathrm{PPV}
(u_{2} ) ) }{n_{D} r_{D,1} r_{D,2} },\vspace*{-2pt}
\end{eqnarray*}
when $u_{1} \leq u_{2}$ and
\begin{eqnarray*}
&& \operatorname{Cov} [ \widehat{\mathrm{PPV}}_{r_{D, 1}, r_{\bar{D},1}}(u_{1}) ,\widehat
{\mathrm{PPV}}_{r_{D, 2}, r_{\bar{D},2}}(u_{2}) ] \\[-3pt]
&&\qquad = \frac{ ( 1- \rho)^2 u_{2} } { (1 - u_{2} )
}\frac{ f_{\bar{D}} (F^{-1} (u_{1} ) ) } {f
(F^{-1} (u_{1} ) ) } \frac{ f_{\bar{D}} (F^{-1}
(u_{2} ) ) } {f (F^{-1} (u_{2} )
) }\\[-3pt]
&&\qquad\quad{}\times \frac{ (r_{D,1} \wedge r_{D,2} ) ( 1 - \mathrm{NPV}
(u_{2}) ) \operatorname{PPV}(u_{1} ) }{n_{D} r_{D,1} r_{D,2} }\\[-3pt]
&&\qquad\quad{} + \frac{ \rho^2 u_{2} } { (1 - u_{2} ) }\frac{ f_{D}
(F^{-1} (u_{1} ) ) } {f (F^{-1} (u_{1}
) ) } \frac{ f_{D} (F^{-1} (u_{2} ) ) } {f
(F^{-1} (u_{2} ) ) } \\[-3pt]
&&\hspace*{11pt}\qquad\quad{}\times \frac{ (r_{\bar{D},2} \wedge
r_{\bar{D},2} ) \operatorname{NPV}(u_{2}) ( 1- \mathrm{PPV} (u_{1} ) )
}{n_{D} r_{D,1} r_{D,2} },
\end{eqnarray*}
when $u_{2} \leq u_{1}$, where $ \sigma^{2}_{\widehat{\mathrm{PPV}}_{r_{D,1},
r_{\bar{D},1}}(u_{1})}$ is defined as in Corollary \ref{PPVNPVuseq}.
\end{corollary}
\begin{pf}
%%%%%%%
Immediate from Theorem \ref{PPVuseqproc}.
%%%%%%%
\end{pf}
\end{appendix}

\section*{Acknowledgments}

The authors would like to thank Jon Wellner for his assistance with the
technical aspects of sequential empirical process theory and Julian
Wolfson for his thoughtful comments on the manuscript. In addition, we
would also like to thank the Associate Editor and referee for their
helpful comments that greatly improved the manuscript.

%suskaldyti doi

% imsref loaded by lrinkeviciute, 2011-12-19 16:10:36
%

\printaddresses

\end{document}